\documentclass{article}
\usepackage{graphicx,amssymb,amsmath,amsthm}

\newtheorem{thm}{Theorem}
\newtheorem{prop}[thm]{Proposition}

\newtheorem{lem}[thm]{Lemma}
\theoremstyle{definition}
\newtheorem{exmp}[thm]{Example}
\newtheorem{defn}[thm]{Definition}
\newtheorem{prob}[thm]{{Questions}}

\title{Tropical Plane Geometric Constructions: a Transfer Technique in Tropical
Geometry}

\author{Luis Felipe Tabera \footnote{The author has been supported by the project MTM2005-08690-C02-02 and a FPU research grant from the Spanish Ministerio de Educaci\'on y Ciencia.}}

\begin{document}
\maketitle
\begin{abstract}
The notion of geometric construction is introduced. This notion allows to
compare incidence configurations in the algebraic and tropical plane. We provide
an algorithm such that, given a tropical instance of a geometric construction,
it computes sufficient conditions to have an algebraic counterpart related by tropicalization. We also provide sufficient conditions in a geometric construction to ensure that the algebraic counterpart always exists. Geometric constructions are applied to transfer classical theorems to the tropical framework, we provide a notion of incidence  theorems and prove several tropical versions of classical theorems like converse Pascal, Fano plane or Cayley-Bacharach.
\end{abstract}

\noindent\emph{keywords:} Tropical geometry, geometric constructions, incidence configurations, classical incidence theorems.

\section{Introduction}
Let $\mathbb{K}$ be an algebraically closed field provided with a non trivial
rank one valuation $v:\mathbb{K}^*\longrightarrow \mathbb{T}\subseteq
\mathbb{R}$, where $v$ is onto $\mathbb{T}$. We have naturally the following map
in the algebraic torus:
\[\begin{matrix}
T: &(\mathbb{K}^*)^n& \longrightarrow& \mathbb{T}^n\\
 &(x_1,\ldots, x_n)& \mapsto& (-v(x_1),\ldots, -v(x_n))
\end{matrix}\]
This map is the \emph{tropicalization} or \emph{projection map}. \emph{Tropical
varieties} are then defined as the image of an algebraic variety $V\subseteq
(\mathbb{K}^*)^n$ under the tropicalization map $T$. Tropical varieties are
polyhedral complexes which are the basic objects of study of tropical geometry.
One of the most interesting aspects of tropical varieties is that they inherit
relevant geometric properties from their algebraic counterparts. In the present
work we try to measure the differences in the behavior of the tropical varieties
with respect to the algebraic ones. In particular, we check if incidence
theorems of classical projective geometry hold in the tropical context. The
origin of this work is the Pappus theorem counterexample showed in \cite{RGST}.
In that paper, it is showed a tropical configuration of points and lines in the
shape of Pappus theorem hypothesis such that it does not verify Pappus thesis.
In particular, it implies that this configuration is not the projection of a
similar configuration of points and lines in the algebraic plane. The authors
provide then another alternative version of the same theorem and claimed that
this new version would hold in the tropical context. The key
of this new version of Pappus theorem is that the hypothesis is
given as the result of a geometric construction dealing with points and lines.
The correctness of this theorem was showed in \cite{Pappus-trop} using some
precursor techniques on geometric constructions. Many incidence theorems can be
given as a construction of a configuration of curves and points
(\emph{hypothesis}) and then some information is derived (the \emph{thesis} of
the theorem). So we will focus on geometric constructions in the plane and how
they behave with respect to tropicalization.

Intuitively, a geometric construction is a procedure that starts with a set of
input curves and points and then define other curves and points by either
intersecting two available curves or computing a curve defined by a polynomial
of fixed support passing through a set of points (a conic through five points,
for example). The main algorithm we present consists in: taking a tropical
instance of a geometric construction, computing a constructible set
$\mathfrak{S}$ over the residual field of the valuation that encodes sufficient
conditions for the compatibility of an algebraic geometric construction. We will also show some certificates during the computation to detect if a tropical realization of a geometric constructions is not the projection of any algebraic realization.

Then, we present the notion of \emph{admissible geometric construction}. This is
a combinatorial notion defined on an associated graph of a geometric construction that
ensure that for all tropical realization, the computed set $\mathfrak{S}$ is non
empty and dense. That is, we will always be able to compute an algebraic
preimage under the tropicalization $T$. This notion can be applied to prove that
some incidence theorems hold in the tropical context if we are able to describe
their hypothesis as the output of an admissible geometric construction. We
provide a notion of \emph{constructible incidence theorem} that is compatible
with tropicalization. In particular, we will show some theorems of this kind.

The notation and basic results are the following: $k$ denotes the
residual field of $\mathbb{K}$ by the valuation. After possibly rescaling the
valuation, we suppose that $\mathbb{Q}\subseteq\mathbb{T}\subseteq\mathbb{R}$.
We will also suppose that we have fixed a multiplicative subgroup $G\subseteq
\mathbb{K}^*$ such that $v: G\rightarrow \mathbb{T}$ is an isomorphism.
$t^\gamma$ denotes the unique element of $G$ such that $v(t^\gamma)=\gamma$. By
the isomorphism, we have that $t^ut^v=t^{u+v}$, $t^0=1$, $t^{-u}=(t^{u})^{-1}$.
$\pi$ denotes the projection from the valuation ring of $\mathbb{K}$ onto $k$.
Let $x\in \mathbb{K}^*$, $u=v(x)$, then $xt^{-u}$ is an element of valuation $0$, so
it has a non zero image in the residual field $k$. We write
$Pc(x)= \pi(xt^{-u})= y\in k^*$ the \emph{principal coefficient} of $x$ to this residual image. Note that the
principal coefficient depends on the group $G$ chosen. The \emph{principal term}
of $x$ is denoted by $Pt(x)=yt^{u}$. This principal term is only a notation, it
is not an element of $\mathbb{K}$ nor $k$ and it cannot be, for example, if the
valuation is a $p$-adic one. It happens that $v(z)=v(x)<v(x-z)$ if and only if
$Pt(z)=Pt(x)$. Usually we will write $x=yt^u+\ldots \quad \textrm{or}\quad
yt^u+o(t^{u})$ in order to emphasize the principal term of an element $x$.

The valuation group $\mathbb{T}$ is given a structure of idempotent semifield
with the tropical operations $``a+b"=\max\{a,b\}$, $``ab"=a+b$. A
tropical polynomial is just a formal sum of monomials $f=``\sum_{i\in
I}a_ix^i"=\max\{a_i+ix: i\in I\}$ where $x=(x_1,\ldots, x_n)$, $i=(i_1,\ldots,
i_n)$, $ix=i_1x_1+\ldots, i_nx_n$. The tropical hypersurface defined by a polynomial is defined as:

\begin{defn}
Let $f=``\sum_{i\in I}a_ix^i" \in\mathbb{T}[x_1,\ldots, x_n]$ be a tropical
polynomial. Then the hypersurface defined by $f$ is the set of points $p\in
\mathbb{T}^n$ such that the value $f(p)=\max\{a_i+ip: i\in I\}$ is attained for
at least two different indices $i,j \in I$. \[\mathcal{T}(f)=\{p: \exists i\neq
j\in I\ \forall k\in I\ a_i+ip=a_j+jp\geq a_k+kp\}\]
\end{defn}

Hypersurfaces defined like this coincide with the projection of algebraic
hypersurfaces by $T$ \cite{einsiedler-2004-}. If $\widetilde{f}=\sum_{i\in
I}\widetilde{a}_ix^i\in \mathbb{K}[x]$ and $f=``\sum_{i\in I} T(a_i)x^i\in
\mathbb{T}[x]"$, then $T(\{\widetilde{f}=0\})=\mathcal{T}(f)$. This result can be
refined
by using residual polynomials.

\begin{defn}
Let $\widetilde{f}=\sum_{i\in I} \widetilde{a}_ix^i \in\mathbb{K}[x]$ be a polynomial in $n$
variables $x=x_1,\ldots,x_n$; $i=i_1,\ldots,i_n$, $Pc(\widetilde{a}_i)= \alpha_i$,
$T(\widetilde{a}_i)=a_i$, $f(x)=``\sum_{i \in I} a_ix^i"$. Let $b=(b_1,\ldots,b_n)\in
\mathbb{T}^n$ be a tropical point. Let
\[\widetilde{f}_{b}(x_1,\ldots,x_n)=\sum_{\substack{i\in I\\ a_i+i_1 b_1+ \cdots + i_n
b_n= f(b_1,\ldots,b_n)} } \alpha_ix^i=Pc(\widetilde{f}(x_1t^{-b_1},\ldots, x_n
t^{-b_n}))\] be the \emph{residual polynomial over $b$}. This is a non zero
polynomial in $k[x_1, \ldots, x_n]$.
\end{defn}

Given a tropical object $x$ (a point, a curve, a configuration,\ldots), a \emph{lift} or \emph{preimage} of $x$ is an algebraic element $\widetilde{x}$ over $\mathbb{K}$ such that $T(\widetilde{x})=x$. In particular, the notion of residual polynomial allows us to compute the lift of a point $x$ belonging to a tropical hypersurface $H$ to an algebraic point $\widetilde{x}$ belonging to a lift $\widetilde{H}$ of $H$. For a constructive proof of this theorem we refer to \cite{Kapranov-EACA} or \cite{Lifting-Constr}.

\begin{thm}\label{initialserie}
Let $\widetilde{f}\in \mathbb{K}[x_1, \ldots, x_n ]$ and
$(\widetilde{b}_1,\ldots, \widetilde{b}_n)\in (\mathbb{K}^*)^n$ be any point,
then there is a root $(\widetilde{c}_1,\ldots, \widetilde{c}_n)$ of
$\widetilde{f}$ such that $Pt(\widetilde{c}_i)=Pt(\widetilde{b}_i)$, $1\leq
i\leq n$, if and only if $b=T(\widetilde{b})$ is a zero of the tropical
polynomial $f$ and $(Pc(\widetilde{b}_1),\ldots,
Pc(\widetilde{b}_n))$ is a root of $\widetilde{f}_b$ in $(k^*)^n$.
\end{thm}

Given a tropical point $q=(q_1,\ldots, q_n)\in \mathbb{T}^n$, it can be written in projective coordinates $q=[q_1:\ldots:q_n:0$. Two tuples $[a_1:\ldots:a_{n+1}]$, $[b_1:\ldots: b_{n+1}]\in \mathbb{T}^{n+1}$ are identified if and only if there is a $c\in \mathbb{T}$ such that $a_i=``cb_i"=c+b_i$, $1\leq i\leq n+1$. Given a tropical polynomial $f=``\sum_{i\in I}a_ix^i"$, it defines
a regular subdivision on its
Newton Polygon that is combinatorially dual to the hypersurface
$\mathcal{T}(f)$. Let $\Delta'$ be the convex hull of the set $\{(i,t)| i\in I,
t\leq a_i\}\subseteq \mathbb{R}^{n+1}$. The upper convex hull of $\Delta'$, that
is, the set of boundary maximal cells whose outgoing normal vector has its last
coordinate positive, projects onto $\Delta$ by deleting the last coordinate.
This projection defines the regular subdivision on $\Delta$.
It is called the \emph{subdivision of $\Delta$ associated to $f$} (See
\cite{Mik05} for the details).

\begin{prop}\label{Newton-Dual}
The subdivision of $\Delta$ associated to $f$ is dual to the set of zeros of
$f$. There is a bijection between the cells of $Subdiv(\Delta)$ and the cells of
$\mathcal{T}(f)$ such that:
\begin{itemize}
\item Every $k$-dimensional cell $\Lambda$ of $\Delta$ corresponds to a cell
$V^\Lambda$ of $\mathcal{T}(f)$ of dimension $n-k$ such that the affine linear
space generated by $V^\Lambda$ is orthogonal to $\Lambda$. (In the case where
$k=0$, the corresponding dual cell is a connected component of
$\mathbb{R}^n\setminus \overline{\mathcal{T}(f)}$)
\item If $\Lambda_1\neq \Lambda_2$, then $V^{\Lambda_1}\cap V^{\Lambda_2} =
\emptyset$
\item If $\Lambda_1\subset \overline{\Lambda}_2$, then $V^{\Lambda_2}\subset
\overline{V^{\Lambda_1}}$
\item $\displaystyle{\mathcal{T}(f)=\bigcup_{0\neq\dim(\Lambda)}V^\Lambda}$
where the union is disjoint.
\item $V^\Lambda$ is not bounded if and only if $\Lambda\subseteq \partial
\Delta$.
\end{itemize}
\end{prop}

In our results, we have to specify families of curves (lines, conics \ldots). A
first approach could be fixing the Newton polygon of the curves. This has a good
geometric meaning. However, without further effort in the proofs, we can fix the
support of the curves. This has no geometric advantages but is a refinement from
an algebraic point of view.

\begin{defn}\label{def_support}
A \emph{support} is a finite subset of $\mathbb{Z}^{n}$ modulo a translation by an
integer vector in $\mathbb{Z}^{n}$. That is, let $\mathcal{P}^f(\mathbb{Z}^n)$ be the
set of finite subsets of $\mathbb{Z}^n$ and let $\sim$ be the relation $A \sim B$ if
and only if there is an integer vector $v\in \mathbb{Z}^n$ such that $A=v+B$. Then,
the set of supports $S(\mathbb{Z}^n)$ of $\mathbb{Z}^n$ is the set of equivalence
classes
$^{\mathcal{P}^f(\mathbb{Z}^n)}\!\diagup\!\sim$. Given a support $I\subseteq
\mathbb{Z}^{n}$, $\delta= \delta(I)$
denotes the number of elements of $I$. $\Delta=cv(I)$, the convex hull of $I$ in
$\mathbb{R}^n$, is the \emph{Newton polytope} of $I$. Note that $\delta$ is
invariant by translations, so it is well defined and $\Delta$ is well defined up
to integer translations. If $H$ is a hypersurface defined by a polynomial
$f=\sum_{i\in I}a_i x^i$, the support of $H$ is the set of tuples $i\in
\mathbb{Z}^n$ such that $a_i$ effectively appear in $f$ modulo integer
translations.
\end{defn}

It is known that in the tropical context, polynomials of different support may
define the same hypersurface. So, we will always fix a priori the support of a
defining polynomial. Sometimes this is not even enough. Contrary to the
algebraic torus, there may be polynomials $f$, $g$ defining the same tropical
curve $C$ but such that one is not a multiple of the other. That is, there is
no monomial $a$ such that $f=ag$. However, for some proofs, it is convenient to
have a tropical polynomial of fixed support defining a curve that is
\emph{canonical} in a sense. We use the notion of concave polynomial from
\cite{Mik05}.

\begin{defn}\label{canonic_polynomial}
To a given tropical polynomial $f=``\sum_{i\in I}a_i x^i"$, we may associate the
function $\varphi: I\subseteq \mathbb{Z}^n \rightarrow \mathbb{T}$, given by
$\varphi(i)=a_i$. We say that $\varphi$ is concave if for any (possibly non
distinct) $i_0, \ldots, i_n\in I\subseteq \mathbb{Z}^n$ and any $t_0,\ldots,
t_n\geq 0$ with $\sum_{k=0}^n t_k=1$ and $\sum_{k=0}^n t_ki_k\in I$ we have
that $\varphi\Big(\sum_{k=0}^n t_ki_k\Big)\geq \sum_{k=0}^{n}t_k\varphi(i_k)$.
If this is the case, we say that $f$ is a \textit{concave polynomial}.
\end{defn}

Fixed the support $I$ and a tropical hypersurface $\mathcal{V}$ defined by a
polynomial $g$ of support $I$, there is (up to a multiplication by a Laurent
monomial) a unique concave tropical polynomial $f$ of support $I$ such that
$\mathcal{T}(f)=\mathcal{V}$.

Finally, we need a notion that is essential in the geometric construction
definition, the notion of stability. Given two tropical curves $C_1$, $C_2$
defined by polynomials $f_1$, $f_2$, it may happen that the intersection of the
curves is infinite even when they share no common component. However, there is
always a finite set of points, called the \emph{stable intersection} of $C_1$
and $C_2$ such that it varies continuously as we deform the coefficients of
$f_1$, $f_2$. This intersection set verifies Bernstein-Kushnirenko theorem,
\cite{RGST}. Namely, let $C_1, C_2$ be two tropical curves defined by
polynomials $f_1$, $f_2$. Let $\Delta_1$, $\Delta_2$ be the Newton polygons of
the respective polynomials. Denote by $\mathcal{M}(\Delta_1, \Delta_2)$ the
mixed volume of $\Delta_1$ and $\Delta_2$, then the number of stable
intersection points of $C_1,C_2$, counted with multiplicities equals
$\mathcal{M}(\Delta_1, \Delta_2)$.

Analogously, given a support $I$ and a set $P$ of $\delta(I)-1$ tropical points,
it may happen that there are infinitely many curves of support $I$ passing
through $P$. However, there is always a unique well defined curve of support $I$
that passes through $P$ and such that it varies continuously as the
configuration $P$ is perturbed. This curve can be computed using tropical linear algebra, see \ref{stable_curve_subsection} and it is called the \emph{stable curve of support $I$ passing through $P$}.

The paper is structured as follows: in
Section~\ref{section_construcciones} we present the notion of
geometric construction and show how to deal with the steps of a construction.
We provide also a notion of \emph{points in
general position} inside a tropical curve that will be useful to generalize some
results. In Section~\ref{section_lift} we provide the main algorithm of
the paper that computes a set of residual sufficient conditions to ensure a
correspondence between algebraic and tropical instances of a geometric
construction. Then, it is shown the limits of the geometric constructions method by a
series of examples, we provide some slight generalizations of the notion of
admissibility and we provide some certificates for the incompatibility with tropicalization.
As well as an example of a construction such that our method cannot derive
neither the compatibility nor the incompatibility because the residual
information is not enough in this case. Finally, in Section~\ref{section_theorem} we use the results obtained so far to provide a notion of incidence theorem that is compatible with tropicalization and to show some instances of theorems of this kind.

\section{The Notion of Geometric Construction}\label{section_construcciones}

We take the notion of incidence structure from the classical context in the study
of finite geometries \cite{Dembowski}. Intuitively, an incidence structure is a
set of points, a set of lines and a set of incidence relations of type
\emph{point $p$ belongs to line $L$}. In our context, we are not only dealing
with lines, but with arbitrary curves in the plane. Still we will control which
curves are accepted in an incidence structure by specifying their support.

\begin{defn}\label{inc_estruct}
A \emph{finite incidence structure} is a tuple $G= (\mathfrak{p}, \mathfrak{B},
\mathfrak{I},Sup)$, where \[\mathfrak{p}\cap \mathfrak{B}=\emptyset, \quad\
\mathfrak{I}\subseteq \mathfrak{p}\times \mathfrak{B}\]
\[Sup:\mathfrak{B}\rightarrow S(\mathbb{Z}^2)\] The elements of $\mathfrak{p}$
are
called \emph{points}, the elements of $\mathfrak{B}$ are \emph{blocks} or
\emph{curves} and the elements of $\mathfrak{I}$ are \emph{flags} or
\emph{incidence relations}. If $x\in \mathfrak{B}$, $Sup(x)\in S(\mathbb{Z}^n)$
is
the support of $x$.
\end{defn}

Every incidence structure $G=(\mathfrak{p}, \mathfrak{B}, \mathfrak{I}, Sup)$ is
naturally identifiable with a labeled graph, the \emph{Levi graph} of the
incidence structure. This is the bipartite graph whose vertices are the elements
of $\mathfrak{p} \cup \mathfrak{B}$ and its edges are the elements of
$\mathfrak{I}$. Each element $x\in \mathfrak{B}$ has as label $Sup(x)$. These
two notions of incidence structures will be used indistinctly.

\begin{figure}[t]
\begin{center}
\includegraphics[width=0.4\linewidth]
{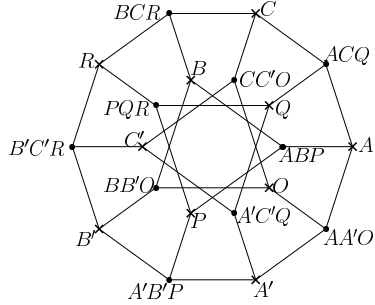}
\caption{The graph of Desargues
configuration}\label{fig:Desargues-grafo-incidencia}
\end{center}
\end{figure}

\begin{exmp}
Desargues Theorem states that two triangles are in perspective with respect to a
point if and only if they are perspective with respect to a line. Desargues
configuration consists in ten points and ten lines. Its incidence structure is:
\[\mathfrak{p}=\{A, B, C, A', B', C', P, Q, R, O\},\]
\[\mathfrak{B}= \{AA'O, BB'O, CC'O, ABP, A'B'P, ACQ, A'C'Q, BCR, B'C'R, PQR\},\]
\[\mathfrak{I}=\{(X_1,X_1X_2X_3), (X_2,X_1X_2X_3), (X_3,X_1X_2X_3)\ |\
X_1X_2X_3\in \mathfrak{B}\}.\]
As every curve in the structure is a line, the support map is constant $Sup(\mathfrak{B})=\{(0,0), (1,0), (0,1)\}$. Figure~\ref{fig:Desargues-grafo-incidencia} represents the incidence graph $G$ of Desargues configuration.
\end{exmp}

\begin{defn}
Let $G=(\mathfrak{p}, \mathfrak{B}, \mathfrak{I}, Sup)$ be an incidence
structure. Denote by $n_\mathfrak{p}$, $n_\mathfrak{B}$ the cardinality of
$\mathfrak{p}$, $\mathfrak{B}$ respectively. For each $y\in \mathfrak{B}$, let
$\delta_y= \delta(Sup(y))$ be the cardinal of the associated support. The
algebraic
support of $G$ is the space \[S_G= \prod_{x \in \mathfrak{p}} (\mathbb{K}^*)^2
\times \prod_{y \in \mathfrak{B}} (\mathbb{K}^*)^{\delta_y-1}.\] The tropical
support of $G$ is the space \[S^t_G= \prod_{x \in \mathfrak{p}} \mathbb{T}^2 \times
\prod_{y \in \mathfrak{B}} \mathbb{T}^{\delta_y-1}.\] We identify the space
$(\mathbb{K}^*)^{\delta_y-1}$ (resp. $\mathbb{T}^{\delta_y-1}$) with the space of
algebraic curves (resp. tropical curves) of support $Sup(y)$ (dehomogenizing the
equation of the curve by a monomial). The dimension of $S_G$ is
$2n_\mathfrak{p}+\sum_{y\in \mathfrak{B}} (\delta_y-1)$.

An algebraic realization (resp. tropical realization) of $G$ is a point
\[(x_1, \ldots, x_{n_p}; y_1, \ldots, y_{n_\mathfrak{B}})\in S_G\ (S_G^t)\]
such that, for every edge $(x_i,y_j) \in \mathfrak{I}$ we have that $x_i\in
y_j$, identifying $y_j$ with the plane  curve (resp. tropical curve) it represents. The
set of algebraic realizations of $G$ is an algebraic set $R_G$ of $S_G$ (resp.
$R_G^t\subseteq S_G^t$).
\end{defn}

A first problem we face at this level is that, in general, $T(R_G)\neq R_G^t$.
This yields the following questions.
\begin{itemize}
\item When does $T(R_G)$ equal $R_G^t$?
\item Given, $x\in R_G^t$, determine if $x$ belongs to $T(R_G)$. In the
affirmative case, compute a preimage $\widetilde{x}$ in $R_G$.
\end{itemize}
In particular, we try to answer these questions using the graph structure of
$G$. This question could be approached using the notion of tropical basis. It
would
consist in taking the equations defining the variety $R_G$. A tropical basis can
be computed from these defining equations (cf. \cite{Computing_trop_var}), the
projection of this basis is a set defining $T(R_G)$, so it would only rest to
check if this basis defines $R_G^t$ or not. This approach together with the
algorithms in \cite{Lifting-Constr} would answer the questions. The main
disadvantage is that they are unfeasible if the graph becomes larger. Here we
will use the strength of the geometric information that the graph contains.

An alternative is to use the graph structure of $G$ and, sometimes, we will not
work with the hole variety $R^t_G$, but with a meaningful subset of it. This
restriction in the set $R^t_G$ is meaningful in the context of geometric
constructions. For the moment, we can derive few information from the graph
structure alone.

\begin{thm}\label{versionincidente_aciclico}
Let $G$ be an incidence structure such that its associated graph is acyclic.
Then, $T(R_G)=R_G^t$. That is, for every tropical realization $x$ of $G$, we can
compute an algebraic realization $\widetilde{x}$ of $G$ that projects correctly
$T(\widetilde{x})=x$.
\end{thm}
\begin{proof}
Let $G$ be the acyclic incidence graph. Reasoning on each connected component of
$G$, we suppose, without loss of generality, that $G$ is a tree. Let $x_0$ be
any node of $G$ and let $\widetilde{x}_0$ be any lift of $x$ to the algebraic
context. The rest of the nodes can be inductively lifted from this one. Let $y$
be an adjacent node to a node $x$ that has already been lifted to
$\widetilde{x}$. We distinguish two cases:
\begin{itemize}
\item $x\in\mathfrak{B}$ and $y\in\mathfrak{p}$. In this case $x$ is a tropical
curve, $\widetilde{x}$ is an algebraic curve projecting onto $x$ and $y$ is a
point in $x$. These are the conditions of Theorem~\ref{initialserie}. Thus,
starting from $y$ we can compute a point $\widetilde{y}$ belonging to
$\widetilde{x}$ and projecting onto $y$.
\item $x$ is a point and $y$ is a curve of support $I=Sup(y)$. $y$ is a tropical
curve of equation $``\sum_{i\in I}a_iz^i"$, with variables $z=(z_1, z_2)$. The
point $\widetilde{x}$ defines, in the configuration space of $\widetilde{y}$,
the hypersurface $H_x$ of curves of support $I$ containing $\widetilde{x}$. Its
equation is $\sum_{i\in I}a_i\widetilde{x}^i$, where the unknowns are the
variables $a_i$. Moreover $y$ belongs to the tropicalization of $H_x$. Thus,
again by Theorem~\ref{initialserie}, it can be computed a lift
$\widetilde{y}$ of $y$ passing through $\widetilde{x}$.
\end{itemize}
\end{proof}

With this Theorem we present a partial answer to the question proposed. However,
acyclic graphs are rather unattractive, because they cannot model many common
situations. Even they cannot deal with the intersection of two conics, because
there will be four intersection points (counted with multiplicities) connected
to both curves and, hence, a cycle in $G$.

If $q$ is a point in a configuration $G$ that belongs to two different curves
$C_1$, $C_2$, it is natural to define $q$ as an intersection point of $C_1$ and
$C_2$. This approach leads to the notion of geometric construction. A
geometric construction is an abstract procedure that produces realizations
(either tropical or algebraic) of an incidence configuration together with an
orientation on the associated graph. Hence, we recall some notation for oriented
(directed) graphs.

A directed graph is a graph such that each edge $\{x_1, x_2\}$ has a defined
orientation $(x_1, x_2)=x_1\rightarrow x_2$. Double orientations in the edges
$x_1\rightarrow x_2$ and $x_2\rightarrow x_1$ are not allowed. For an oriented
edge $x_1\rightarrow x_2$, we say that $x_1$ is a \emph{direct predecessor} of
$x_2$ and that $x_2$ is a \emph{direct successor} of $x_1$. An \emph{oriented
path} is a chain of oriented edges $x_1 \rightarrow x_2 \rightarrow \ldots
\rightarrow x_n$. If there is an oriented path from $x_1$ to $x_n$, we say that
$x_1$ is a \emph{predecessor} of $x_n$ and that $x_n$ is a \emph{successor} of
$x_1$. An \emph{oriented cycle} is an oriented path such that its starting node
equals its ending node, $x_1=x_n$. A directed graph without oriented cycles is
called a \emph{directed acyclic graph} (DAG). If $G$ is a DAG, the nodes $x$ of
$G$ that are not the successor of any other node are called \emph{sources}. Any
node $x$ of a DAG $G$ has associated a \emph{depth}. If $x$ is a source then its
depth is $0$. If $x$ is not a source, let $y_1,\ldots, y_n$ be the direct
predecessors of $x$. The depth of $x$ is defined as:
$\textrm{depth}(x) = 1 + \max\{\textrm{depth}(y_1), \ldots, \textrm{depth}(y_n)\}.$ The depth of a DAG $G$ is the maximal depth of its nodes.

\begin{defn}\label{construccion}
A geometric construction is an abstract procedure consisting in:
\begin{itemize}
\item Input elements: two finite subsets $\mathfrak{p}_0$, $\mathfrak{B}_0$ such
that $\mathfrak{p}_0\cap \mathfrak{B}_0= \emptyset$ and a support map
$Sup: \mathfrak{B}_0 \rightarrow S(\mathbb{Z}^2)$. Initially, the set of incidence relations is the empty set $\mathfrak{I}=\emptyset$.
\item Steps of the construction, a finite sequence of different steps:
\begin{itemize}
\item Given a support $I$ with $\delta(I)=n\geq 2$ and $n-1$ points $\{q_1,
\ldots, q_{n-1}\}$ we add a new curve $C$ of support $I$ to $\mathfrak{B}$, we
also add new oriented incidence conditions $q_i\rightarrow C$, $1\leq i\leq
n-1$.
\item Given two curves $C_1$, $C_2$ of support $I_1$, $I_2$ and Newton Polygons
$\Delta_1$, $\Delta_2$ respectively, we add
$M=\mathcal{M}(\Delta(I_1),\Delta(I_2))$ new source points $q_1,\ldots,q_M$. We
add the oriented incidence conditions $C_1\rightarrow q_i$, $C_2 \rightarrow
q_i$, $1 \leq i\leq M$.
\end{itemize}
\item Output: an incidence graph $G$ provided with an orientation.
\end{itemize}

A \emph{tropical realization of a geometric construction} $\mathfrak{C}$ is a
tropical realization of its associated graph $G$ such that:
\begin{itemize}
\item If $x\in \mathfrak{B}$ is a curve and it is not an input element, let $I$
be its support and let $\{y_1,\ldots, y_{\delta(I)-1}\}$ be the direct
predecessors of $x$. Then $x$ is exactly the stable curve of support $I$ passing
through the set of points $\{y_1,\ldots,y_{\delta(I)-1}\}$.
\item If $x\in \mathfrak{p}$ and it is not an input point, let $y_1$, $y_2$ be
the direct predecessors of $x$ and let $\{x_1,\ldots, x_n\}$ be the common
direct successors of $y_1$ and $y_2$. Then, $\{x_1, \ldots, x_n\}$ are exactly
the stable intersection of $y_1$ and $y_2$, counted with multiplicities.
\end{itemize}

An \emph{algebraic realization of a geometric construction} $\mathfrak{C}$ is an
algebraic realization of its associated graph $G$ such that:
\begin{itemize}
\item If $x\in \mathfrak{B} \setminus \mathfrak{B}_0$, let $I$ be its support
and let $\{y_1,\ldots, y_{\delta(I)-1}\}$ be the direct predecessors of $x$.
Then, $x$ is the unique algebraic curve of support $I$ that passes through the
points $\{y_1,\ldots, y_{\delta(I)-1}\}.$
\item If $x\in \mathfrak{p}$ and it is not an input point, let $y_1$, $y_2$ be
the direct predecessors of $x$ and let $\{x_1,\ldots, x_n\}$,
$n=\mathcal{M}(\Delta_1,\Delta_2)$ be the common direct successor of $y_1$ and
$y_2$. Then, the curves $y_1$, $y_2$ intersect exactly in the finite set of
points $\{x_1, \ldots, x_n\}$ where the points are counted with multiplicities.
\end{itemize}
\end{defn}

Given an algebraic (resp. tropical) realization of the input elements of a
geometric construction $\mathfrak{C}$, there can only be finitely many
realizations of $\mathfrak{C}$ with these input elements, because the
realizations of the rest of the elements are fixed by the input elements and the
steps of the construction. The only possibility to have different realizations
of $\mathfrak{C}$ with the same input elements is a permutation of the labels of
the intersection (resp. stable intersection) of two curves $y_1$, $y_2$ and the
consequent changes in the successor elements of $y_1, y_2$ in the construction.

It is clear that, in the tropical plane, every step of a construction can be
performed. That is, given two curves $C_1$, $C_2$, we can always define the set
of $\mathcal{M}(\Delta_1, \Delta_2)$ intersection points (counted with
multiplicities). Analogously, the stable curve through a set
of points is always well defined. Thus, in the tropical context, given a
tropical realization of the input elements of $\mathfrak{C}$, there is always a
realization of $\mathfrak{C}$ with these input elements. However, this is not
the case in the algebraic context. Two different curves $C_1$, $C_2$ may share a
common component. Here, we cannot define a finite intersection set with the nice
properties the tropical stable intersection has. Even if the intersection set of
the curves is finite, there may not be enough intersection points in the torus.
For example, the lines $3x+2y+4$, $5x+y+2$ do not have any intersection point in
the torus. These degenerate cases should be avoided. So, we need a notion of a
well defined construction. A geometric construction is \emph{well defined} if it
is well defined for a generic realization of the input elements. That is, let
$R_{0}$ be the space of algebraic realizations of the input elements
$\mathfrak{p}_0\cup \mathfrak{B}_0$. In this case, as the set of incidence
conditions is empty, the realization space equals the support space,
$R_{0}=S_{0}$. Let $L$ be the set of configurations such that every step of the
construction $\mathfrak{C}$ is well defined (that is, the projection into $R_0$
of the algebraic realizations of $\mathfrak{C}$). The construction $G$ is
\emph{well defined} if $L$ is dense in $R_{0}$.

It is clear that the oriented graph $G$ of a geometric construction
$\mathfrak{C}$ never has an oriented cycle, so $G$ is always a directed acyclic
graph (DAG). The input elements are exactly the sources and every node of $G$
has defined a depth. Usually, proofs are made by induction on the depth of $G$.

In practice, many interesting incidence configurations can be defined as a
subgraph of the graph of a geometric construction. Sometimes we will have to add
additional elements to fit the incidence configuration into the definition of
geometric construction. Hence, we present a characterization of the incidence
graphs $G$ that appear as a subgraph of a geometric construction.

\begin{prop}\label{subgrafodeconstruccion}
Let $G$ be an incidence graph provided with an orientation. Then it is the
subgraph of the graph of a geometric construction if and only if
\begin{itemize}
\item $G$ is a directed acyclic graph, (\emph{DAG}).
\item If $x$ is a vertex of type $\mathfrak{p}$, then it has at most two direct
predecessor.
\item If $x$ is a curve of support $I$, then $x$ has at most $\delta(I)-1$
direct predecessors.
\item If $x,y$ are two different curves with a common direct successor, then
they have at
most $\mathcal{M}(\Delta_x,\Delta_y)$ common direct successors.
\item If $x$ and $y$ are two curves with the same support $I$ and both curves
have exactly $\delta(I)$ direct predecessor, then the sets of direct
predecessors are different.
\end{itemize}
Moreover, $G$ is exactly the graph of a geometric construction if and only if
the previous inequalities are equalities for every node different from a source.
\end{prop}
\begin{proof}
Let $G$ be a graph satisfying all these conditions, a construction
$\mathfrak{C}$ can be defined such that it contains $G$ as a subgraph. Every
source of $G$ is defined as an input element. Suppose defined the construction
of every element of depth up to $i$, the definition of the depth $i+1$ elements
is as follows. Let $x$ be a point ($x\in \mathfrak{p}$) of depth $i+1$, if it has two predecessors
$y,z$, then they have at most $\mathcal{M}(\Delta_y,\Delta_z)$ common direct
successors. If there are not enough intersection points, we add points of depth
$i+1$ up to $\mathcal{M}(\Delta_x,\Delta_y)$ and define all of them (in
particular $x$) as the intersection of $y$ and $z$. If $x$ is a point of depth
$i+1$ that has only one direct predecessor $y$, we add a line $z$ as an input
curve (a curve of support $\{(0,0), (1,0),(0,1)\}$), define it as a direct
predecessor of
$x$ and proceed as in the previous case. In the case where $x$ is a curve of
support $I$ and depth $i+1$, there are at most $\delta(I)-1$ predecessors of
$x$. Add to the construction $\mathfrak{C}$ as many input points as necessary up
to $\delta(I)-1$ and define $x$ as the curve passing through these points. Note
that the last condition of the hypothesis disallow the construction to have
repeated steps. If two curves $x$ and $y$ of the same support $I$ have both
$\delta(I)$ direct predecessors, then the set of direct predecessors is
different, so $x$ and $y$ are curves obtained by different steps.

This method defines a construction $\mathfrak{C}$ that contains $G$ as a
subgraph. It is clear that $G$ is exactly the graph of $\mathfrak{C}$ if and
only if the equalities in the hypothesis hold.
\end{proof}

One might be tempted to add additional allowed steps to a construction besides
the two steps of the definition. In particular, a common step in
Classical Geometry is to choose a point in a curve.
Proposition~\ref{subgrafodeconstruccion} proves that this step does not increase
the expressivity of the constructions. If $\mathfrak{C}$ is a geometric
construction such that the additional step of taking a curve through a point or
taking a point inside a curve is allowed, then the graph of $\mathfrak{C}$ is
the subgraph of another construction $\mathfrak{C}_1$ without these additional
steps. So, in practice, we may work with this additional step with the agreement
that \emph{``choosing a point in a curve is essentially equivalent to add an
input line (curve of support $\{(0,0),(1,0),(0,1)\}$) to our construction,
intersect the line with the curve and choose one intersection point.''} See for
example Theorem~\ref{pascal_converse} for an example of this technique of adding
additional elements to a familiar incidence configuration in order to obtain a
geometric construction.

The advantage of the construction method over a direct approach to the study of
incidence configurations is that the problem is almost reduced to lifting the
steps of
the construction.

\subsection{The Stable Curve Through a Set of Points}\label{stable_curve_subsection}

Consider now the problem of lifting the curve of support $I$ passing through a
set of points. Either in the algebraic or tropical context, this curve can be
computed solving a linear system of equations. Let $q_1,\ldots, q_{\delta-1}$ be
the set of points we want the curve to pass through. Let $f=``\sum_{i\in
I}a_ix^{i_1}y^{i_2}"$ be a polynomial defining the curve of support $I$ passing
through the set of points. The coordinates $a_i$ of $f$ belong to the
hyperplanes defined by $``\sum_{i\in I}z_i q_{j1}^{i_1}q_{j2}^{i_2}"$, $1\leq
j\leq \delta-1$. Thus, the coordinates $a_i$ form a solution of a homogeneous tropical
linear system of equations. The stable intersection of the hyperplanes can be
computed using tropical Cramer's rule \cite{RGST}. This stable intersection of hyperplanes is exactly the coordinates of the stable curve $f$. In order to lift these linear systems of equations, we recall the following basic facts of tropical linear algebra:

A \emph{tropical matrix} of dimension $n\times m$ is a matrix with coefficients
in $\mathbb{T}$. The \emph{tropical determinant} of a square matrix is defined as:
\[\begin{vmatrix}
x_{11}& \ldots & x_{1n}\\
\vdots& & \vdots\\
x_{n1}& \ldots & x_{nn}
\end{vmatrix}_t=``\sum_{\sigma \in \Sigma_n}x_{1\sigma(1)} \cdots x_{n\sigma(n)}" =
\mathop{max}_{\sigma \in \Sigma_n}\{ x_{1\sigma(1)} + \cdots +
x_{n\sigma(n)}\}\]
where $\Sigma_n$ is the permutation group of $n$ elements. A square tropical matrix
is called \emph{singular} if the value of its tropical determinant is attained
for at least two different permutations $\sigma$ and $\tau$. In other case it is
called \emph{regular}.

Tropical and algebraic determinants can be related by the notion of
pseudodeterminant. Let $A=(a_{ij})$ be a $n\times n$ tropical matrix. Let
$B=(b_{ij})$ be a $n\times n$ matrix with coefficients over any ring $R$. Let
$|A|_t$ be the tropical determinant of $A$. We define:
\[\Delta_A(B)=\sum_{\substack{\sigma\in\Sigma_n\\
``a_{1\sigma(1)} \ldots a_{n,\sigma(n)}" = |A|_t}}
(-1)^{i(\sigma)}b_{1\sigma(1)} \cdots b_{n\sigma(n)}\]
the \emph{pseudodeterminant} of $B$ with respect to weight $A$. With this notion
we can derive sufficient conditions for the compatibility of the algebraic and
tropical determinant.

\begin{defn}\label{pseudocramer}
Let $A=(a_{ij})$ be a $n\times(n\!+\!1)$ tropical matrix. Let $B=(b_{ij})$ be a
matrix with coefficients in a ring $R$ with the same dimension as $A$. We denote
\[\textrm{Cram}_A(B)=(S_1,\ldots,S_{n+1})\]
where $S_i=\Delta_{A^{i}}(B^{i})$ and $A^{i}$ (respectively, $B^{i}$) denotes
the corresponding submatrix obtained by deleting the $i$-th column in $A$
(respectively, $B$).
\end{defn}

\begin{lem}\label{pseudocramerlema} Suppose we are given a system of $n$
linear homogeneous equations in $n\!+\!1$ variables in $\mathbb{T}$. Let $A$ be
the coefficient matrix of the system. Let $\widetilde{A}$ be any matrix with
coefficients in $\mathbb{K}$ such that $T(\widetilde{A})=A$. Let $B=Pc(\widetilde{A})$ be the
matrix of principal coefficients of $\widetilde{A}$. If no element of
$\textrm{Cram}_A(B)$ vanishes, then the linear system defined by $\widetilde{A}$
has only one projective solution and its tropicalization equals the stable
tropical solution $[ | A^{1} |_t : \ldots : |A^{n\!+1\!} |_t ]$.
\end{lem}
\begin{proof}
See \cite{Pappus-trop}
\end{proof}

If one pseudodeterminant $\Delta_{A^i}(B^i)=0$, there is a lack of information
of what the principal coefficient of the determinant $|A^i|$ is and, more
serious, the control on the tropicalization $T(|A^i|)$ is lost. A careful look
at these badly behaved systems yields the following:

\begin{prop}\label{certificado_curva_por_puntos}
Let $A$ be a $n\times n+1$ tropical matrix. Let $x=[\ |A^1|_t : |A^2|_t : \ldots
: |A^{n+1}|_t\ ]$ be the stable solution of the linear system of equations
defined by $A$. Let $\widetilde{A}$ be any matrix in $\mathbb{K^*}$ projecting onto
$A$ and $B=Pc(A)$. Let $\textrm{Cram}_A(B)=(S_1,\ldots, S_{n+1})$. Then:
\begin{itemize}
\item If every tropical determinant $|A^i|_t$ is regular, then $S_i\neq 0$, the
homogeneous linear system defined by $\widetilde{A}$ has only one solution
$\widetilde{x}$ and it projects onto $x$, $T(\widetilde{x})=x$.
\item If $S_j= 0$ and there is an index $i$ such that $S_i\neq 0$, then the
homogeneous linear system $\widetilde{A}$ has only one projective solution
$\widetilde{x}$, that never tropicalizes correctly: $T(\widetilde{x})\neq x$.
\item If $S_i=0$ for all $i$, we do not have any information. The linear system
defined by $\widetilde{A}$ may be either determined or undetermined. If
$\widetilde{x}$ is a solution of the system, both possibilities
$T(\widetilde{x})=x$ and $T(\widetilde{x})\neq x$ can occur, even if the
solution $\widetilde{x}$ is unique.
\end{itemize}
\end{prop}
\begin{proof}
If $A^i$ is regular, then $|A^i|_t=``a_{1,j_1}\cdots a_{n, j_{n}}"$ is attained
for only one permutation. It follows that $\Delta_A(B)= b_{1,j_1}\cdots
b_{n,j_n}\neq 0$ for any matrix $B$ with entries in $k^*$. Hence, the
algebraic system is determined, because at least the $i$-th projective
coefficient $|\widetilde{A}^i|$ is not zero. Moreover, in this case it will
always happen that $T(|\widetilde{A}^i|)=|A^i|_t$. If every tropical matrix
$A^i$ is regular, then we have the first item.

For the second item, if $S_j=0$, then $T(|\widetilde{A}^j|)<|A^j|_t$. It is even
possible that $|\widetilde{A}^j| = 0$. But, as $S_i\neq 0$, then
$T(|\widetilde{A}^i|)=|A^i|_t$, so the coefficient $i$ can be used to
dehomogenize. If follows that $\widetilde{x}$ is well defined (because
$|\widetilde{A}^i|\neq 0$), but it cannot projects into $x$ because they will
always differ in the term $j$.

Finally, in the case where $S_i=0$ for every $S$ we cannot decide if the system
is determined without further information. This depends on the terms of higher
order of the elements of $\widetilde{A}$. For an illustrative example, let
$\mathbb{K}$ be the field of Puiseux series, let
\begin{center}
\begin{tabular}{rlrl}
$A=$&$\begin{pmatrix}0&0&0\\0&0&0\end{pmatrix}$&
$\widetilde{A}_1=$&$\begin{pmatrix}1&1&1\\1&1&1\end{pmatrix}$\\
$\widetilde{A}_2=$&$\begin{pmatrix}1+t&1+t^2&1+t^3\\1&1&1\end{pmatrix} $&
$\widetilde{A}_3=$&$\begin{pmatrix}1+t&1+2t&1+3t\\1&1&1\end{pmatrix}$
\end{tabular}
\end{center}
The three matrices $\widetilde{A}_1$, $\widetilde{A}_2$, $\widetilde{A}_3$
projects into $A$. All of them satisfy that
\[\textrm{Cram}_{A}(Pc(\widetilde{A}))=(0,0,0).\]
The tropical stable solution of the tropical system is the point $[0:0:0]$. The
first algebraic system $\widetilde{A}_1$ is undetermined and it contains points
such that $\widetilde{x}=[1:1:-2]$ that projects correctly onto $[0:0:0]$ and
other points such that $\widetilde{x}=[1:t:-1-t]$ that does not. The second
system $\widetilde{A}_2$ is a determined system such that its unique solution
$\widetilde{x}=[t^2-t^3: -t+t^3: t -t^2]$ does not project into $x$. The last
system $\widetilde{A}_3$ is a determined one. Its solution is $[-1:
2: -1]$ and projects correctly.
\end{proof}

Before establishing the relationship of the algebraic and tropical curve, let us
check some properties of the pseudodeterminants. From
Lemma~\ref{pseudocramerlema}, it follows that if the entries of the matrix $B$
are indeterminates, then no pseudodeterminant $\Delta_A(B)$ vanishes and the
algebraic determinant projects correctly. However, it may happen that the
entries of the matrix $B$ are algebraically dependent elements. For example,
suppose we are computing the conic $a_{xx} x^2+a_{yy} y^2+a_{xy}xy + a_x x+ a_y
y+ a_1$ passing through a set of points $\{p_1,\ldots,p_5\}$, $p_1=(b_1,b_2)$.
This conic can be computed using linear algebra. In the matrix $B$ that
describes the linear system to solve, the terms $b_1^2$,
$b_2^2$, $b_1b_2$ will appear in the system of equations. These monomials are
not algebraically independent. Nevertheless, in order to apply
Lemma~\ref{pseudocramerlema}, it is only needed that the involved
pseudodeterminants do not vanish. Now it is proved that, if the residual
coefficients $(\gamma_1, \gamma_2)$ of the points $p_1$ are indeterminates (or
generic elements), then, the pseudodeterminants are never zero. The next is a
rather technical Lemma that proves a stronger property.

\begin{lem}\label{Levantamientounpaso}
Let $C_i=\{c_i^1,\ldots,c_i^{j_i}\}$, $1 \leq i\leq r$ be disjoint sets of
variables. Suppose that we have $F_u = \{f_u^1,\ldots,f_u^{n+1}\} \subseteq
k[\bigcup_{i=1}^r C_i]$, $1\leq u \leq n$ sets of polynomials in the variables
$c_i^j$. Suppose also that the following properties hold:
\begin{itemize}
\item For a fixed set $F_u$, $f_u^l$, with $1\leq l\leq n+1$ are
multihomogeneous polynomials in the sets of variables
$C_{u^1},\ldots,C_{u^{s_u}}$ with the same multidegree.
\item If $u\neq v$ then $F_u$, $F_v$ involve different sets of variables $C_i$.
\item In a family $F_u$, if $l\neq m$ then the monomials of $f_u^l$ are all
different from the monomials of $f_u^m$.
\end{itemize}
Let us construct the $n\times (n\!+\!1)$ matrix
\[B=(f_u^l)_{\substack{1\leq u\leq n,\\ 1\leq l\leq n+1}}\]
Let $A$ be any $n\times (n\!+\!1)$ tropical matrix. Write
\[S=\textrm{Cram}_A(B)=(S_1,\ldots,S_{n+1}).\] Then
\begin{enumerate}
\item $S_1,\ldots,S_{n+1}$ are non identically zero multihomogeneous polynomials
in the sets of variables $C_1,\ldots,C_r$ with the same multidegree.
\item If $\sigma, \tau$ are different permutations in $\Sigma_{n+1}$ which
appear in the expansion of $S_l$ (and, therefore $\sigma(n+1)=\tau(n+1)=l$),
then all resulting monomials in $\prod_{u=1}^n (A^{l})_u^{\sigma(u)}$ are
different from the monomials in $\prod_{u=1}^n (A^{l})_u^{\tau(u)}$
\item If $l\neq m$, then $S_l$, $S_m$ have no common monomials.
\end{enumerate}
\end{lem}
\begin{proof}
See \cite{Pappus-trop}
\end{proof}

In the case of computing the algebraic curve $\widetilde{C}$ through a set of
points $\widetilde{P}$, suppose for simplicity that the points $\widetilde{q}_i$
are given
in homogeneous coordinates with generic principal coefficients and
tropicalization $[q_i^1 : q_i^2: q_i^3]$. \[\widetilde{q_i}=[\gamma_i^1
t^{-q_i^1}+ \cdots : \gamma_i^2 t^{-q_i^2}+ \cdots: \gamma_i^3 t^{-q_i^3}+
\cdots].\] Suppose also that the defining equation of $\widetilde{C}$ is
homogenized adding a new variable $z$, \[\widetilde{C}\equiv \sum_{i\in I}
\widetilde{a}_i x^{i^1} y^{i^2} z^{r-i^1-i^2}.\] Let $\widetilde{A}$ be the
matrix of this homogenized linear systems and $B=Pc(\widetilde{A})$. We claim
that the matrix $B$ is in the conditions of Lemma~\ref{Levantamientounpaso}. The
$j$-th row of $B$ is \[B_j=\big((\gamma_j^1)^{i_1^1} (\gamma_j^2)^{i_1^2}
(\gamma_j^3)^{r-i_1^1 - i_1^2}, \ldots, (\gamma_j^1)^{i_\delta^1}
(\gamma_j^2)^{i_\delta^2} (\gamma_j^3)^{r-i_\delta^1- i_\delta^2}\big)\] Hence,
in the hypothesis of Lemma~\ref{Levantamientounpaso}, $C_i=\{\gamma_j^1,
\gamma_j^2, \gamma_j^3\}$, each polynomial $f_u^l$ is a different homogeneous
monomial. So, the hypothesis holds. Thus, we conclude that for this homogenized
system, the vector $\textrm{Cram}_A(B)$, that contains a representative of the
residues of the vector of coefficients of $\widetilde{C}$, belongs to the torus,
$\textrm{Cram}_A(B)\in (k^*)^n$. It follows that the algebraic solution
$[\widetilde{a}_{i_1}: \ldots: \widetilde{a}_{i_\delta}]\in \mathbb{K}^*$. Finally,
as every coefficient of every point $\widetilde{q}_j$ and $[\widetilde{a}_{i_1}:
\ldots: \widetilde{a}_{i_\delta}]$ is nonzero, we can dehomogenize everything.
The pseudodeterminants $\Delta_{A^i}(Pc(\widetilde{A}^i))$ are nonzero provided
that $Pc(\widetilde{q}_i)= (\gamma_i^1,\gamma_i^2)$ are generic. To sum up, we
have the following:

\begin{thm}\label{condiciones_de_pseudocramer}
Let $I$ be a support, $\delta=\delta(I)$, $P=\{q_1,\ldots, q_{\delta-1}\}$,
$q_j=(q_j^1, q_j^2)$ a set of tropical points, $\widetilde{P}=\{
\widetilde{q}_1, \ldots, \widetilde{q}_{\delta -1}\}$ a set of algebraic points
such that $\widetilde{q}_j=(\widetilde{q}_j^1, \widetilde{q}_j^2) =(\gamma_j^1
t^{- q_j^1} + \ldots, \gamma_j^2 t^{- q_j^2} + \ldots)$. Let $C$ be the stable
tropical curve of support $I$ passing through $P$ computed using Cramer's rule.
Let
\[A=((q_j^1)^{i^1}(q_j^2)^{i^2})\quad \quad
\widetilde{A}=((\widetilde{q}_j^1)^{i^1}(\widetilde{q}_j^2)^{i^2})\]
be the matrices of the linear system defining $C$ and $\widetilde{C}$. For
simplicity, it is assumed that the columns of $A$ are indexed by the set $I$.
Then, the pseudodeterminants are non identically zero polynomials in the set
$\{\gamma_j^i, 1\leq j\leq \delta-1, 1\leq i\leq 2\}$. If the pseudodeterminants
verify that
\[\Delta_{A^i}(Pc(\widetilde{A}^i))\neq 0,\quad i\in I.\]
then, there is only one curve $\widetilde{C}$ passing through $\widetilde{P}$
and $T(\widetilde{C})=C$. That is, the pseudodeterminants provide residual
sufficient conditions for the equality $T(\widetilde{C})=C$.

In this case, let $\widetilde{f}=\sum_{i\in I} \widetilde{a}_i x^{i^1}y^{i^2}$
be the polynomial of support $I$ defining $\widetilde{C}$ computed by Cramer's
rule, suppose that this polynomial is dehomogenized with respect to the index
$i_0$ ($\widetilde{a}_{i_0}=1$), then, the principal coefficients of
$\widetilde{a}_i$ are
\[(Pc(\widetilde{a}_1), \ldots, Pc(\widetilde{a}_\delta))=\Big( \frac{
\Delta_{A^{i_1} }
(Pc(\widetilde{A}^{i_1 }))}{ \Delta_{A^{i_0 }} (Pc(\widetilde{A}^{i_0 }))},
\ldots,
\frac{\Delta_{A^{i_\delta} }(Pc(\widetilde{A}^{i_\delta }))}{ \Delta_{A^{i_0 }}
(Pc(\widetilde{A}^{i_0 }))} \Big)\]

\end{thm}
\begin{proof}
If no pseudodeterminant $\Delta_{A^i}(Pc(\widetilde{A}^i))$ vanishes, then
$T(|\widetilde{A}^i|)=|A^i|_t$. In particular, no determinant
$|\widetilde{A}^i|$ is zero. Let
\[\widetilde{C}\equiv \{\sum_{i\in I}|\widetilde{A}^i|x^{i_1}y^{i_2}=0\}\]
be the unique algebraic curve of support $I$ passing through $\widetilde{P}$ and
projecting onto $C$, the curve defined by $``\sum_{i\in I}|A^i|_t
x^{i_1}y^{i_2}"$, i.e. the stable tropical curve through $P$.

Note that if no pseudodeterminant vanishes, the coordinates of $\widetilde{C}$
belongs to the algebraic torus in homogeneous coordinates $(\mathbb{PK}^*)^\delta$.
Thus, if one wants an affine representation of the coordinates of the curve, it
can be dehomogenized with respect to any index $i_0\in I$ and still the result
will project correctly into the (dehomogenized) equation of the tropical curve
$C$. Furthermore, taking principal coefficients commutes with dehomogenization
in $(\mathbb{PK}^*)^\delta$, so the last claim holds.
\end{proof}

We have shown sufficient conditions for the compatibility of the algebraic and
tropical curve through a set of corresponding points. If the lifts of points
$\widetilde{P}$ are residually generic, the algebraic curve $\widetilde{C}$
passing through them is unique. We know that this curve projects onto the stable
curve through the tropical points, but it is not clear what is the residual
relationship of its coefficients.

This is important in the context of incidence configurations. Proofs such as the
one in Theorem~\ref{versionincidente_aciclico} are done recursively in the graph
of the configuration. So, if using residually generic coefficients is an
argument to Theorems such as \ref{condiciones_de_pseudocramer} and we want to
use this Theorem in an induction scheme, we should establish the residual
genericity of the coefficients of the curve $\widetilde{C}$. Next, we prove that
if the points $\widetilde{q}_i$ are residually generic, then the coefficients of
$\widetilde{C}$ are also residually generic.

\begin{thm}\label{la_curva_por_puntos_genericos_es_generica}
Let $I=\{l_1, \ldots, l_\delta\}$, $l_k=(i_k,j_k)$ be a support. Let $P= \{q_1,
\ldots, q_{\delta-1}\}$ be a set of tropical points. Let $C$ be the stable
tropical curve of support $I$ passing through $P$. Let
$\widetilde{P}=\{\widetilde{q}_1, \ldots, \widetilde{q}_{\delta-1}\}$,
$Pc(\widetilde{q}_i)=(\gamma_i^1, \gamma_i^2)$ and $\widetilde{C}$ the algebraic
curve of support $I$ passing through $\widetilde{P}$. Let
$\widetilde{f}=\sum_{(i,j)\in I}\widetilde{a}_{i,j}x^{i}y^{j}$ be the algebraic
curve representing $\widetilde{C}$ dehomogenized with respect to the index
$l_0=(i_0,j_0)$. Let $\gamma_1=\{\gamma_1^1, \ldots, \gamma_{\delta-1}^1\}$,
$\gamma_2=\{\gamma_1^2, \ldots, \gamma_{\delta-1}^2\}$ Consider the map
\[\begin{matrix}
k^{2\delta -2}&\longrightarrow & k^{\delta -1}&\\
(\gamma_1,\gamma_2)&\mapsto&
\textrm{Cramer}(\gamma_1,\gamma_2)=&\Big( \frac{ \Delta_{A^{l_1} }
(Pc(\widetilde{A}^{l_1 }))}{ \Delta_{A^{l_0 }} (Pc(\widetilde{A}^{l_0 }))},
\ldots,
\frac{\Delta_{A^{l_\delta} }(Pc(\widetilde{A}^{l_\delta }))}{ \Delta_{A^{l_0 }}
(Pc(\widetilde{A}^{l_0 }))} \Big)
\end{matrix}\]
that represents the principal coefficients of $\widetilde{f}$ in terms of the
principal coefficients of $\widetilde{P}$ (provided that the pseudodeterminants
of Theorem~\ref{condiciones_de_pseudocramer} do not vanish). Then, the map
$\textrm{Cramer}$ is dominant, that is, if the principal coefficients of
$\widetilde{P}$ are generic, then the polynomial $\widetilde{f}$ is generic
among the polynomials of support $I$ dehomogenized with respect to $l_0$.
\end{thm}

\begin{proof}
Write $q_l=(q_l^{1},q_l^{2})$, $C=\mathcal{T}(``\sum_{ij}a_{ij}x^iy^j")$. Then,
$C$ is the curve defined by the stable solution of:
\[``\sum_{(i,j)\in I}a_{ij}(q_l^{1})^i(q^{2}_l)^j", 1\leq l\leq \delta-1\]
and the lifts of $C$ verify the relations
\[\sum_{(i,j)\in
I}\widetilde{a}_{ij}(\widetilde{q}_{l}^1)^i(\widetilde{q}_{l}^2)^j=0, 1\leq
l\leq \delta-1\]
Take the equations
\[\widetilde{f}_l=\sum_{(i,j)\in
I}\widetilde{a}_{ij}x^iy^jt^{-a_{ij}-iq^{1}_l-jq^{2}_l}, 1\leq
l\leq \delta-1,\]
which correspond to a (tropical) translation of the problem to the point $0$. We
dehomogenize the tropical equation of $C$ $(a_{i_0j_0}=0)$, and the algebraic
equation of $\widetilde{C}$ ( $\widetilde{a}_{i_0j_0}=1$) with respect to a
term $(i_0, j_0)\in I$. The residual conditions on the principal coefficients
$\alpha_{ij}$ of $\widetilde{a}_{ij}$ are:
\[f_l=\sum_{J_l}\alpha_{ij}(\gamma_{l}^1)^i(\gamma_{l}^2)^j, 1\leq l\leq
\delta-1,\]
where $J_l\subseteq I$ are the monomials such that $-a_{ij} -iq^{1}_l -jq^{2}_l$
is minimized. Notice that, by construction, each $J_l$ has at least two terms.
Write $\alpha=\{\alpha_{ij}| (i,j)\neq (i_0,j_0)\}$, $\gamma_1= \{\gamma_{1}^1,
\ldots, \gamma_{\delta-1}^{1}\}$, $\gamma_2 =\{\gamma_{1}^2, \ldots,
\gamma_{\delta-1}^{2}\}$. Each residual equation $f_l$ is affine in the set of
variables $\alpha$, and the coefficients of this affine equations are monomials
in $\{\gamma_{l}^{1},\gamma_{l}^{2}\}$. Moreover, we know that there are nonzero
solutions to this system. Without loss of generality, every polynomial $f_l$ can
be saturated with respect to the coordinate hyperplanes (that is, we eliminate
redundant $\gamma$). These polynomials are still denoted by $f_l$. Thus, we have
a system of equations in $3\delta -3 $ unknowns.

Let $\mathcal{V}$ be the Zariski closure of the image of the map:
\[\begin{matrix}\label{param_variedad_cramer}
k^{2\delta -2}&\longrightarrow & k^{3\delta -3}\\
(\gamma_1,\gamma_2)&\mapsto&
(\gamma_1,\gamma_2,\textrm{Cramer}(\gamma_1,\gamma_2))
\end{matrix}\]
It is clear that this is a birational map between the space $k^{2\delta-2}$ and
$\mathcal{V}$. Let $\mathcal{I}$ be the ideal of $\mathcal{V}$. $\mathcal{I}$ is
a prime ideal that contains the polynomials $(f_1, \ldots, f_{\delta -1})$ in
$k[\alpha,\gamma_1,\gamma_2]$. By construction, the field of rational functions
of $\mathcal{V}$ is the field of fractions of the integer domain
\[\mathbb{L}=\textrm{Frac}\left(\frac{k[\gamma_1,\gamma_2,\alpha]}{\mathcal{I}}
\right)=k(\gamma_1,\gamma_2)\]
In particular, $\gamma_1,\gamma_2$ is a transcendence basis of
$k\subseteq\mathbb{L}$ and the dimension of $\mathbb{L}$ is $2\delta-2$. For
each $f_l$, if the variable $\gamma_{l}^1$ does not appear in $f_l$, then
$\gamma_{l}^2$ is an element of $\mathbb{L}$ which is algebraic over
$k(\alpha,\gamma_{l}^1)$. Analogously, if $\gamma_{l}^2$ does not appear in
$f_l$, then $\gamma_{l}^1$ is algebraic over $k(\alpha, \gamma_{l}^2)$. If both
variables appear in $f_l$, then just choose $\gamma^{j}_l$ algebraic over
$k(\alpha,\gamma^{3-j}_l)$. In this way, the set $g=\alpha\cup
\{\gamma^{3-j}_l$, $1\leq l\leq \delta -1\}$ is such that $\mathbb{L}$ is
algebraic over $k(g)$. As $\#g=2\delta -2$, we conclude that $g$ is a
transcendence basis of $k\subseteq \mathbb{L}$. In particular, the set $\alpha$
is algebraically independent over $k$. This means that:
\begin{equation}\label{corte_con_los_ejes}
\mathcal{I}\cap k[\alpha]=\mathcal{I}\cap k[\gamma_1,\gamma_2]=0
\end{equation}
Hence, the projection of $\mathcal{V}$ over the space of coordinates $\alpha$ is
dense in $k^{\delta-1}$. But the image of the projection is the image of
$k^{2\delta-2}$ by the map $\textrm{Cramer}$, so $\textrm{Cramer}$ is dominant.
\end{proof}

\subsection{Points in Generic Position in a Curve}\label{section_points_general}
Before dealing with the problem of the intersection of two curves, let us
explore the notion of points in general position inside a curve. This notion
will be helpful for subsequent results and is directly related with the notion
of stable curve through a set of points. First, an adequate notion of tropical
points in general position must be provided. There are slightly different
approaches to this definition in the literature. All of them share the same
idea, but apply to different problems, see for example \cite{Mik05},
\cite{Tesis-Hannah}, or \cite{GM-Numbers-of-curves}. These notions are
adequate for enumerative problems, but not for the incidence structures we
study. Moreover, we want to provide a notion of generic points in a fixed curve
$C$. Informally, a set of points $P$ is in general position inside a curve $C$
if $C$ is the unique curve of its type that contains $P$. Again, to formalize
this we use the notion of stability:

\begin{defn}\label{mi_nocion_de_puntos_en_posicion_general}
Let $C$ be a tropical curve of support $I$. A set of points $q_1, \ldots, q_n$,
$n\leq \delta(I)-1$ is in \emph{generic position with respect to
$C$} if there are tropical points
$q_{n+1}, \ldots, q_{\delta-1}$ such that $C$ is the stable curve of support $I$
passing through $q_1,\ldots, q_{\delta-1}$.
\end{defn}

One would like to characterize the points in general position in a curve $C$
because, in general, it is not easy to check the Definition. A first result is
the following:

\begin{lem}
Let $C$ be a curve of support $I= \mathbb{Z}^2 \cap \Delta$, where $\Delta$ is a
convex polygon. Suppose that the dual subdivision induced by $C$ in $\Delta$ is
a triangulation that has all points in $\Delta \cap I$ as vertices. Let $q_1,
\ldots, q_{\delta-1}$ be different points in $C$ such that every point $q_i$
lies in the relative interior of an edge of $C$ and two different points do not
lie in the same edge. Let $\Gamma$ be the graph contained in the subdivision of
$\Delta$ consisting of those edges such that their dual contains a point $q_i$.
If $\Gamma$ is a maximal tree contained in $Subdiv(\Delta)$, then the vertices of
$\Gamma$ are exactly the points of $I$ and $C$ is the unique curve of support
$I$ passing through $q_1, \ldots, q_{\delta-1}$. In particular, $q_1, \ldots,
q_{\delta-1}$ are points in general position in $C$
\end{lem}
\begin{proof}
We refer to \cite{Mik05}.
\end{proof}

This Lemma only works for very special curves, because of the restriction on the
support of the curve and the induced subdivision (triangulation) in $\Delta$.

\begin{defn}
Let $C$ be a tropical curve of support $I$ and Newton Polygon $\Delta$. Let
$\Gamma_0$ be the skeleton of $Subdiv(\Delta)$ associated to $C$ (the set of
cells of dimension $0$ and $1$. This is always a connected graph). We modify
$\Gamma_0$ adding to $\Gamma_0$ every point in $I\setminus \Gamma_0$
as follows.

If $x_1, \ldots, x_r \in I$ are the points of $I$ lying in the interior of an
edge
$e$ of $\Gamma_0$, then we add these points as $2$-valent vertices of $\Gamma_0$
splitting the edge $e$ into $r+1$ edges. If $x\in I$ lies in the relative
interior of a polygon $\Delta_v$ of the subdivision, then $x$ is added to
$\Gamma_0$ as an isolated point. In this case, the resulting graph, denoted by
$\Gamma$, is no longer connected.

Let $q$ be a point in $C$. If $q$ lies in an edge of $C$, let $\Delta_q$ be the
dual edge in $\Gamma_0$, then $\Delta_q=e_1\cup \ldots \cup e_d$ is refined as a
union of edges in $\Gamma$. An \emph{assignment} of $q$ is a
choice of one of
the edges $e_1, \ldots, e_d$. In the case where $q$ is a vertex of $C$, the dual
cell $\Delta_q$ of this vertex is a polygon. Let $S$ be the set of isolated
points of $I$ in the interior of $\Delta_q$ and $e_1, \ldots, e_d$ be the set of
refined edges in the boundary $\Gamma\cap \partial \Delta_q$. An
\emph{assignment} of $q$ is a choice of an element in $S \cup \{e_1, \ldots,
e_d\}$.

If $q_1,\ldots, q_n$ are points (possibly repeated) in $C$, an \emph{assignment
of the points} is an assignment of each point $q_i$ such that:
\begin{itemize}
\item Let $q_{i_1}, \ldots q_{i_r}$ be the points lying in the same edge of $C$,
let $\Delta_{q}=e_1\cup \ldots \cup e_d$ be the refined dual edge in $\Gamma$.
It is required that the assignment of $q_{i_j}$ is different from the assignment
of $q_{i_k}$ whenever $j\neq k$ (even in the case that $q_{i_j}=q_{i_k}$ is a
repeated point).
\item Let $q_{i_1}, \ldots, q_{i_r}$ be points identified with a vertex (that
is, a vertex with multiplicity $r$). Let $\Delta_q$ be the polygon dual to the
vertex. Let $l=\# \{ \Delta_q \cap I \}$. It is required that at most $l$ points
are assigned to different points in $S$ and that the $r-l$ other points are
mapped to different refined edges of the boundary of $\Delta_q$ that has not been previously assigned.
\item The set of refined edges of $\Gamma$ such that have assigned a point $q_i$
form an acyclic subgraph $\Lambda$ of $\Gamma$.
\end{itemize}
\end{defn}

\begin{lem}\label{alcanzan_puntos}
Let $C$ be a curve of support $I$. Let $q_1, \ldots, q_{\delta-1}$ be a list of
points such that there exists an assignment in $\Gamma$. Then
\begin{itemize}
\item Every point of $I$ that lies in the relative interior of a polygon
$\Delta_v$ of $Subdiv(\Delta)$ is assigned to a point $q_i$.
\item The set of assigned edges $\Lambda$ is a maximal tree in $\Gamma$ that
contains as
vertices every non isolated vertex of $\Gamma$.
\end{itemize}
\end{lem}
\begin{proof}
The proof is based on the properties of lattice subdivisions of tropical
curves presented in \cite{Mik05}. Let $S$ be the set of points of $I$ lying in
the relative interior of a polygon in $Subdiv(\Delta)$ and let $l$ be the number
of these points. Let $r=\delta-l$ be the number of non isolated vertices of
$\Gamma$. Then, at most $l$ points $q_i$ are assigned to a point in $S$ and at
least $\delta-1-l=r-1$ points are assigned to an edge on $\Gamma$. Then, from
the property that the set of assigned edges of $\Gamma$ is an acyclic graph. It
follows that the number of assigned edges must be smaller than the number of
vertices. That is, the number of assigned edges must be exactly $r-1$. It
follows that the graph of assigned edges is connected, i.e. a tree. Moreover,
this tree is maximal, because it attains every non isolated vertex of $\Gamma$.
Finally, the number of isolated points of $\Gamma$ assigned to a point is $l$
(every isolated point has been assigned).
\end{proof}

\begin{lem}\label{No_hace_falta_asignar_vertices}
Let $C$ be a tropical curve of support $I$ and Newton polygon $\Delta$. Let
$\Gamma$ be the refinement of $\Gamma_0$. Let $q_1,\ldots,q_{\delta-1}$ be
points in the curve. Suppose that if a vertex $v$ of $C$ coincides with $r$
points $q_i$, then the dual polygon $\Delta_v$ contains exactly $r$ point of $I$
in its interior. Suppose that there is an assignment of the points. Then, $C$ is
the stable curve passing through $q_1,\ldots,q_{\delta-1}$.
\end{lem}
\begin{proof}
Let $\widetilde{q}_i$ be lifts of the points $q_i$ with generic residual
coefficients $\gamma_j=(\gamma_{j}^{1},\gamma_{j}^{2})$. In order to define a
curve $\widetilde{C}$, we have to compute lifts of the coefficients
$\widetilde{a}_i$ of a polynomial defining $C$. Let $f$ be the concave
polynomial of support $I$ defining $C$, $f=``\sum_{i\in I}a_ix^{i^1}y^{i^2}"$
dehomogenized with respect to a vertex $i_0$ of the polygon $\Delta$
($a_{i_0}=0$). Notice that, if $g=``\sum_{i\in I}b_ix^{i^1}y^{i^2}"$ is any
tropical polynomial of support $I$ such that if $i$ is a vertex of $Subdiv(C)$
then $b_i=a_i$ and $b_i\leq a_i$ in any other case, then $f$ and $g$ represents
the same piecewise affine function and $\mathcal{T}(g)=C$. We will compute a
polynomial $g$ with this characteristic.

Given an edge $e$ of $Subdiv(\Delta)$, let $e=e_1 \cup\ldots \cup e_{d-1}$ be
the refinement in $\Gamma$, $e_k=[i_k,i_{k+1}]$. If there were two different
edges $e_k, e_l$, $k<l$ that are not assigned to any point $q_j$, then, if
$k+1=l$ then the vertex $i_{k+1}$ would be a vertex of $\Gamma$ that is not
attained by $\Lambda$, if $k+1<l$ then either $\Lambda$ does not attain a vertex
of $\Gamma$ (if $e_{k+1}, \ldots, e_l$ are not assigned) or $\Lambda$ is not
connected (if at least one $e_j$ is assigned with $k<j<l$), contrary to the
results in Lemma~\ref{alcanzan_puntos}. Hence, for the case of an edge
$\Delta_q=e_1 \cup \ldots \cup e_{d}$, at most one of the refined edges $e_k$ is
not assigned to any point. The residual values $\alpha_i$ for a point $i$ of $I$
contained in an edge of $Subdiv(\Delta)$ are computed recursively, starting from
$\alpha_{i_0}=1$. By the maximal tree structure of $\Lambda$ as a subgraph of
$\Gamma$, we can always
suppose that we are in one of the following two cases:

\textit{1})The edge is $e=[{i_1}, \ldots {i_d}]$, we only know the value of
$\alpha_{i_1}$ and there are exactly $d-1$ points $q_{j_1}, \ldots q_{j_{d-1}}$
in the
dual edge $V^e\subseteq C$. The non homogeneous residual system of equations
associated to the points is:
\[\left \{ \begin{matrix}
\alpha_{i_1}\gamma_{j_1}^{i_1}&+&\cdots&+&\alpha_{i_d}\gamma_{j_1}^{i_d}&=0\\
\alpha_{i_1}\gamma_{j_2}^{i_1}&+&\cdots&+&\alpha_{i_d}\gamma_{j_2}^{i_d}&=0\\
&\multicolumn{3}{c}{\cdots\cdots}\\
\alpha_{i_1}\gamma_{j_{d-1}}^{i_1}&+&\cdots&+&\alpha_{i_d}\gamma_{j_{d-1}}^{i_d}
&=0
\end{matrix}\right.\]
in the unknowns $\{\alpha_{i_2}, \ldots, \alpha_{i_d}\}$ and
$\gamma_{j_l}^{i_l}=(\gamma_{j_l}^1)^{i_l^1}(\gamma_{j_l}^2)^{i_l^2}$. This
system is determined. To show this, we may homogenize each row of the monomial
matrix $(\gamma)$ by a new variable $\gamma_{l_i}^3$, hence, we obtain a matrix
that is in the hypothesis of Lemma~\ref{Levantamientounpaso}. By this Lemma, we
conclude that
its minors are non identically zero multihomogeneous polynomials that will
remain non identically zero after dehomogenizing each variable
$\gamma_{l_i}^3=1$. The determination of $\alpha_{i_1}$ is just a
dehomogenization of the solution. Hence, we conclude that there is only one
solution $\{\alpha_{i_2}, \ldots, \alpha_{i_d}\}$ of this linear system in the
algebraic torus over the residual field $(k^*)^{d-1}.$ Notice that, using
induction, each $\alpha_{i_j}$ is a non zero rational function in $\alpha_{i_0}$
and $\gamma$. Applying this steps recursively we can compute the values of every
edge of integer length $d-1$ and $d-1$ assigned points. Notice that, in
particular, we can compute the values of every $\alpha_i$ associated to a vertex
of $Subdiv(\Delta)$ and that they are non zero.

\textit{2}) The edge is $e=[a_{i_1},\ldots,a_{i_d}]$ and the values of
$\alpha_{i_1}$ and $\alpha_{i_d}$ have been already computed. Necessarily, there
are exactly $d-2$ points $q_{j_1},\ldots, q_{j_{d-2}}$ in the dual edge of $e$,
because if there where more points, there would be a cycle in the graph
$\Lambda$, contrary to the hypothesis, and if there where less points, $\Lambda$
would not be a maximal tree. The residual conditions on the unknowns
$\{\alpha_2, \ldots, \alpha_{d-1}\}$ form a non homogeneous system of $d-2$
linear equations in $d-2$ unknowns with a similar structure to the previous
case. So, if the coefficients of $\gamma_i$ are generic, there is only one
solution (this time in $k^{d-2}$ because the determination of the values of
$\alpha_{i_1}$ and $\alpha_{i_d}$ do not correspond to just a dehomogenization).
Again, applying induction, each $\alpha_i$ is a rational function of
$\alpha_{i_0}$ and $\gamma$.

Thus, if the coefficients $\gamma$ are generic, all the values $\alpha_i$
corresponding to an index $i$ that is not an isolated vertex of $\Gamma$ can be
computed from $\gamma$ and $\alpha_{i_0}$ and its value is unique. It only rest
to compute the values $\alpha_i$ corresponding to indices in $I$ belonging to
the relative interior of a polygon in $Subdiv(\Delta)$. In this case, the
corresponding point $q_i$ lie in a vertex $v \in C$. Let $\Delta_v$ be its dual
polygon in $Subdiv(\Delta)$. Every coefficient corresponding to $\partial
\Delta_v \cap I$ has been already computed. Let $\{{j_1}, \ldots, {j_r}\}=
\partial\Delta_v \cap I$ and $\{{k_1}, \ldots, {k_s}\} = int(\Delta_v)\cap I$.
There are $s$ points $q_i$ identified to $v$. The residual system of equations
corresponding to these points is:
\[\left \{ \begin{matrix}
\alpha_{{k_1}}\gamma_{l_1}^{{k_1}}&+\cdots+&\alpha_{{k_s}}\gamma_{l_1}^{{k_s}}
&=-\alpha_{{j_1}}\gamma_{l_1}^{{j_1}}&-\cdots-&\alpha_{{j_r}}\gamma_{l_1}^{{j_r}
}\\
\alpha_{{k_1}}\gamma_{l_2}^{{k_1}}&+\cdots+&\alpha_{{k_s}}\gamma_{l_2}^{{k_s}}
&=-\alpha_{{j_1}}\gamma_{l_2}^{{j_1}}&-\cdots-&\alpha_{{j_r}}\gamma_{l_2}^{{j_r}
}\\
&\multicolumn{3}{c}{\cdots\cdots}\\
\alpha_{{k_1}}\gamma_{l_s}^{{k_1}}&+\cdots+&\alpha_{{k_s}}\gamma_{l_s}^{{k_s}}
&=-\alpha_{{j_1}}\gamma_{l_s}^{{j_1}}&-\cdots-&\alpha_{{j_r}}\gamma_{l_s}^{{j_r}
}\\
\end{matrix}\right.\]
in the unknowns $\{\alpha_{{k_1}}, \ldots, \alpha_{{k_s}}\}$. Again, if the
values of $\gamma$ are generic, there is only one solution in $k^s$.

So, starting from the value $\alpha_{i_0}=1$ the rest of the values are
determined from $\gamma$. Let $\widetilde{a}_i$ be any element of $\mathbb{K}^*$
such that if $\alpha_i\neq 0$ then $Pt(\widetilde{a}_i)=\alpha_i t^{-a_i}$, and,
if $\alpha_i=0$, then $Pt(\widetilde{a}_i)=t^{-a_i+1}$ (\emph{``a posteriori"}
one could show that this later case does not happen). Let
$\widetilde{g}=\sum_{i\in
I}\widetilde{a}_i x^{i^1}y^{i^2}$. Let $\widetilde{C}$ the algebraic curve
defined by $\widetilde{g}$, its projection $T(\widetilde{C})$ is the curve $C$.
But it may happen that $\widetilde{C}$ does not contain the points
$\widetilde{q}_i$, because the computations have been done just in the residual
field. Anyway, by construction, the principal terms of $\widetilde{q}_i$ are in
the hypothesis of Theorem~\ref{initialserie}, we can compute points
$\widetilde{q}_i'$ lying in $\widetilde{C}$ such that
$Pt(\widetilde{q}_i')=Pt(\widetilde{q}_i)$. That is, there is a curve
$\widetilde{C}$ passing through a set of lifts $\widetilde{q}'_i$ of $q_i$ with
generic residual coefficients in the sense of
Theorem~\ref{condiciones_de_pseudocramer}. Hence, $C=T(\widetilde{C})$ is the
stable curve passing through $q_1,\ldots,q_{\delta-1}$.
\end{proof}

\begin{thm}\label{combinatoria_de_posicion_generica}
Let $C$ be a curve of support $I$ and Newton polygon $\Delta$, let $\Gamma$ be
the refinement of the subdivision $\Delta$. Let $q_1, \ldots, q_{\delta-1}$ be
points in the curve. If there is an assignment of $q_1, \ldots, q_{\delta-1}$,
then $C$ is the stable curve of support $I$ passing through the points.
\end{thm}

\begin{proof}
For each vertex $v$ of $C$ containing points $q_{j_1}, \ldots, q_{j_r}$, let
$q_{j_{s+1}}, \ldots q_{j_r}$ be the points identified with $v$ but that are
assigned to an edge of $\Gamma$ and
let $e_1,\ldots e_r$ be those edges. Perturb the point $q_{j_i}$ in $C$
translating it along the dual edge of $e_i$. Denote this point by $q_{j_i}'$.
For the rest of points, take $q_{j_i}'=q_{j_i}$. The points $q'_1,\ldots
q'_{\delta-1}$ are points in $C$ in the conditions of
Lemma~\ref{No_hace_falta_asignar_vertices}. Hence, $C$ is the stable curve
through $\{q_{1}',\ldots, q'_{\delta-1}\}$. Making a limit process on each
perturbed point $q_{j_i}'\rightarrow q_{j_i}$ along the edge that contains the point, the stable curve $C$ trough the points $\{q_{1}',\ldots, q'_{\delta-1}\}$ stays invariant along the whole process. By the continuity of the stable curve through perturbations of a set of points, we conclude that $C$ is the stable curve through $q_1,\ldots,q_{\delta-1}$.
\end{proof}

It is conjectured that the conditions imposed in the preceding Theorem are also
necessary in order to have the genericity of the points inside the curve. That
is, we claim that given $C$ a tropical curve and $q_1,\ldots, q_{\delta-1}\in
C$, $C$ is the stable curve through the points if and only if there is an
assignment of the points. In many concrete examples it can be easily shown that
this condition is a complete characterization of a set of points in general
position in a curve. But the problem is still open for an arbitrary curve.

\subsection{Stable Intersection of Curves}
Let us face the second kind of steps in a geometric construction, namely the
intersection of two curves given by two polynomials $f$, $g$. In this case we
find a similar result as in the case of the curve passing through a set of
points. Given two lifts $\widetilde{f}$, $\widetilde{g}$, some polynomials in
the residual coefficients of $f$ and $g$ can be computed such that, if none of
them vanish, then the intersection of $\widetilde{f}$ and $\widetilde{g}$ is a finite set of
points projecting onto the stable intersection of $f$ and $g$. To obtain this
result and compute the stable intersection itself we use the notion of tropical
resultant \cite{Resultantes-trop}.

The tropical resultant of two univariate polynomials with fixed support is
defined as the tropicalization of the algebraic resultant of two generic
polynomials of the same support.

\begin{defn}
Let $I$, $J$ be two finite subsets of $\mathbb{N}$ of cardinality at least 2
such that $0\in I\cap J$. That is, the support of two polynomials that do not
have zero as a root. Let $R(I,J,\mathbb{K})$ be the resultant of two polynomials
with indeterminate coefficients, $\widetilde{f}=\sum_{i\in I} a_ix^i$, $\widetilde{g}=\sum_{j\in J}
b_jx^j$ over the field $\mathbb{K}$.
\[R(I,J,\mathbb{K})\in \mathbb{Z}/(p\mathbb{Z})[a,b],\]
(where $p$ is the characteristic of the field $\mathbb{K})$. Let $R_t(I, J,
\mathbb{K})$ be the tropicalization of $R(I,J,\mathbb{K})$. This is a polynomial
in $\mathbb{T}[a,b]$, which is called the tropical resultant of supports $I$ and
$J$ over $\mathbb{K}$.
\end{defn}

For the bivariate case, the tropical resultant is defined as the specialization
of the adequate univariate resultant substituting the variables by
univariate polynomials:

\begin{defn}
Let $\widetilde{f}$ and $\widetilde{g}$ be two bivariate polynomials. In order
to compute the algebraic resultant with respect to $x$, we can rewrite them as
polynomials in $x$.
\[\widetilde{f}=\sum_{i\in I}\widetilde{f}_i(y)x^i, \quad
\widetilde{g}=\sum_{j\in J}\widetilde{g}_j(y)x^j,\] where
\[\widetilde{f}_i=\sum_{k=o_i}^{n_i}A_{ik}t^{-\nu_{ik}}y^k,\
\widetilde{g}_j=\sum_{q=r_j}^{m_j}B_{jq}t^{-\eta_{jq}}y^q\] and $A_{ik}$,
$B_{jq}$ are residually generic elements of valuation zero (indeterminates). Let
$P(a_i,b_j,\mathbb{K})=R(I,J,\mathbb{K})\in \mathbb{Z}/(p\mathbb{Z})[a_i,b_j]$
be the algebraic univariate resultant of supports $I$, $J$. The algebraic
resultant of $\widetilde{f}$ and $\widetilde{g}$ is the polynomial
$P(\widetilde{f}_i, \widetilde{g}_j, \mathbb{K})\in \mathbb{K}[y]$. Analogously,
let $f=T(\widetilde{f})$, $g=T(\widetilde{g})$, $f=``\sum_{i\in I}f_i(y)x^i"$,
$g=``\sum_{j\in J}g_j(y)x^j"$, where \[f_i= ``\sum_{k=o_i}^{n_i} \nu_{ik} y^k",\
g_j= ``\sum_{q=r_j}^{m_j} \eta_{jq} y^q".\] Let $P_t(a_i, b_j, \mathbb{K})
=R_t(I, J, \mathbb{K}) \in \mathbb{T} [a_i, b_j]$ be the tropical resultant of
supports $I$ and $J$. Then, the polynomial $P_t(f_i, g_j, \mathbb{K})\in
\mathbb{T}[y]$ is the \emph{tropical resultant} of $f$ and $g$.
\end{defn}

The tropical resultant polynomials $R_t(I,J,\mathbb{K})$ and
$P_t(f_i,g_j,\mathbb{K})$ depend on the characteristics of the fields
$\mathbb{K}$ and $k$. However, the hypersurfaces they define do not depend on
the characteristics of the fields. Moreover, they provide a method to relate the
stable intersection of two curves. We refer to \cite{Resultantes-trop} for the
details.

Suppose that we are given two bivariate tropical polynomials $f,g\in
\mathbb{T}[x,y]$ of support $I,J$ respectively. Let $R_x(x)$ be the resultant of
$f$ and $g$ with respect to $y$. Then, the tropical roots of $R_x(x)$ are
exactly the $x$ coordinates of the stable intersection of $C_1=\mathcal{T}(f)$
and $C_2=\mathcal{T}(g)$. If we compute also $R_y(y)$ the tropical resultant of
$f$ and $g$ with respect to $x$, we have that $P=\mathcal{T}(R_x(x))\cap
\mathcal{T}(R_y(y))$ is always a finite set containing the stable intersection
of $C_1, C_2$. Let $a$ be a natural number such that $x-ay$ is injective in $P$.
Let $R_z(z)$ be the resultant of the polynomials $f(zy^a,y), g(zy^a,y)$ with
respect to $y$. Then, we have that $C_1\cap C_2\cap \mathcal{T}(R_x(x))\cap
\mathcal{T}(R_y(y))\cap \mathcal{T}(R_z(xy^{-a}))$ is exactly the stable
intersection of $C_1$ and $C_2$. In order to ensure the compatibility of the
stable intersection with the algebraic intersection, we just compute residually
sufficient conditions for the compatibility of the resultants. Let $\widetilde{f}$, $\widetilde{g}$ be two lifts of $f$
and $g$ with residually generic coefficients. Let
$\widetilde{R}_x(x)=Res(\widetilde{f},\widetilde{g},y)$, $\widetilde{R}_y(y)=Res(\widetilde{f},\widetilde{g},x)$, $\widetilde{R}_z(z)=
Res(\widetilde{f}(zy^a, y), \widetilde{g}(zy^a,y), y)$. Then $\widetilde{R}_x(x) =\sum_{i\in K_1}
\widetilde{a}_ix^i$, $R_x(x) =\sum_{i\in K_1} a_ix^i$ and $T(\widetilde{a}_i)\leq a_i$. The
tropical polynomial $R_x$ induces a subdivision in the convex hull of $K_1$,
which is an interval in $\mathbb{R}$ with integer endpoints. We have that
$T(V(\widetilde{R}_x(x))) =\mathcal{T}(R_x(x))$ if, for every index $j$ corresponding to
a vertex of the subdivision induced by $R_x(x)$, $T(\widetilde{a}_j)=a_j$, which is equivalent to say that
the principal coefficient of $\widetilde{a}_j$ is $\alpha_j t^{-a_j}$. But $\alpha_j$ is
a polynomial in the residual coefficients of $\widetilde{f}, \widetilde{g}$. Let
$\{\alpha_j\}\cup \{\beta_j\} \cup \{\gamma_j\}$ be the polynomials in the
residual coefficients of $\widetilde{R}_x(x)$, $\widetilde{R}_y(y)$, $\widetilde{R}_z(z)$ corresponding to vertices of the subdivision of their Newton polytopes. If none of them vanish, there
will be a correspondence between the algebraic and tropical resultant. Moreover,
this provides a relation between the algebraic intersection and the tropical
stable intersection. In particular:

\begin{thm}\label{interestable}
Let $\widetilde{f}$, $\widetilde{g}\in\mathbb{K}[x,y]$. Then, it can be computed
a finite set of polynomials in the residual coefficients of $\widetilde{f}$,
$\widetilde{g}$ depending only on their tropicalization $f$, $g$ such that, if
no one of them vanish, the tropicalization of the intersection of
$\widetilde{f}$, $\widetilde{g}$ is exactly the stable intersection of $f$ and
$g$. Moreover, the multiplicities are conserved.
\[\sum_{\substack{\widetilde{q}\in \widetilde{f}\cap \widetilde{g}\\
T(\widetilde{q})=q}} \textrm{mult} (\widetilde{q}) =\textrm{mult}_t(q)\]
\end{thm}

So the step of intersecting two curves is also compatible with tropicalization
in the residually generic case. The problem we face now in order to use this
result in a nontrivial geometric construction is to determine the residual
genericity of the intersection points of the two curves. Of course, it is not
true in general that the intersection points of two curves are points in general
position. A classical example is the intersection set $P$ of two generic cubics
in the plane. In this case, $P$ has $9$ points and all of them lie on two
different cubics. As there is only one cubic passing through $9$ points in
general position, it follows that $P$ cannot be a set of points in general
position. Actually eight of the points determine the ninth, \cite{MR1376653}.
However, taking strict subsets of $P$, it is expected that these sets of points
are in general position. This is the aspect we want to explore. The election of
adequate subsets of the intersection points is done by geometric properties
of the corresponding tropical intersection points.

\begin{thm}\label{teorema_de_puntos_interseccion_en_posicion_general}
Let $C_1$, $C_2$ be two curves of support $I_1$, $I_2$ and Newton polytopes
$\Delta_1$, $\Delta_2$ respectively. Let $q=\{q_1, \ldots, q_n\}$ be a set of
points contained in the stable intersection of $C_1$ and $C_2$ such that $q$ is
in general position (Definition~\ref{mi_nocion_de_puntos_en_posicion_general})
with respect to both curves. Let $\widetilde{C}_1$, (respectively $\widetilde{C}_2$) be a lift
of $C_1$ (resp. $C_2$), expressed by a polynomial $\widetilde{f}$, (resp. $\widetilde{g}$) of
support $I_1$, (resp. $I_2$) and dehomogenized with respect to an index $i_0$,
(resp. $j_0$) that is a vertex of the Newton Polygon $\Delta_1$, (resp.
$\Delta_2$). Suppose that the residual coefficients of the polynomials $\widetilde{f}$, $\widetilde{g}$ range over a dense Zariski open subset of $k^{\delta_1+ \delta_2 -2}$ and let $\widetilde{q}_i$ be lifts of the points $q_i$ to the intersection of the algebraic curves. Then, the tuple of possible values of $(Pc(\widetilde{q}_1), \ldots, Pc(\widetilde{q}_n))$ contains an open dense subset of $k^{2n-2}$. That is, if the residual coefficients of $\widetilde{f}$ and $\widetilde{g}$ are generic, so they are the tuple of coefficients of $\widetilde{q}_i$.
\end{thm}
\begin{proof}
Let
\[f_1 =``\sum_{(i_1, i_2) \in I_1} a_i x^{i_1} y^{i_2}" \hspace{1cm}
f_2= ``\sum_{ (j_1, j_2) \in I_2} b_{j} x^{j_1} y^{j_2}"\]
be two tropical polynomials defining $C_1$ and $C_2$ and let
\[\widetilde{f}_1= \sum_{(i_1, i_2) \in I_1} \widetilde{a}_i x^{i_1} y^{i_2} \hspace{15mm}
\widetilde{f}_2 =\sum_{ (j_1, j_2) \in I_2} \widetilde{b}_j x^{j_1} y^{j_2}\]
be the lifts of the curves. Without loss of generality, it is supposed that both
polynomials are dehomogenized with respect to two monomials that are vertices of
$\Delta_1$ and $\Delta_2$ respectively.
Let $\alpha_i= Pc(\widetilde{a}_i)$, $\beta_j= Pc(\widetilde{b}_j)$, $(\gamma_{1l},
\gamma_{2l})= Pc(\widetilde{q}_l)$, $\alpha=\{\alpha_i\}$, $\beta=\{\beta_j\}$,
$\gamma=\{\gamma_{kl}\}$. As the points are in general position, it must be the
case $n\leq \min\{\delta_1,\delta_2\} -1 $. The proof mimics the reasoning of
Theorem~\ref{la_curva_por_puntos_genericos_es_generica}. So, a parametrization
of the residual coefficients of the curves and the points $\widetilde{q}_i$ is needed.
The local equations $(\widetilde{f}_1)_{q_i}$, $(\widetilde{f}_2)_{q_i}$ form a linear system
of equations in the residual coefficients of the points $\gamma_{il}$ where the unknowns are
the residual coefficients of the curves $\alpha_i$, $\beta_j$. This is a linear system of $2n$
equations in at most $\delta_1+\delta_2-2$ unknowns of full rank. It follows
that we may take $\alpha_0=\{\alpha_{i_1}, \ldots, \alpha_{i_{\delta_1-n-1}}\}$
residual coefficients of $\widetilde{f}_1$ as parameters such that the remaining system
is determined. Analogously, we may take $\beta_0=\{\beta_{i_1}, \ldots,
\beta_{i_{\delta_2-n-1}}\}$ residual coefficients such that the remaining system
of equations in determined. It follows that the remaining variables $\alpha_i$,
$\beta_j$ are rational functions of $\alpha_0$, $\beta_0$ and $\gamma$. These
rational functions define the parametrization
\[\begin{matrix}
k^{\delta_1+\delta_2-2}& \rightarrow & k^{\delta_1+\delta_2+2n-2}\\
(\alpha_0, \beta_0, \gamma)& \mapsto& (\alpha, \beta, \gamma)
\end{matrix}\]
of a variety $\mathcal{V}$ that can be identified with the vectors of principal
coefficients $(C_1, C_2, q)$. Let $\mathbb{L}$ be the field of fractions of
$\mathcal{V}$. It is clear that every class $\gamma_{ki}$ is algebraic over
$k(\alpha,\beta)\subseteq \mathbb{L}$ and that $\mathbb{L}=k(\alpha_0,\beta_0,\gamma)$
by the parametrization. Thus, $\{\alpha_0,\beta_0,\gamma\}$ and $\{\alpha,
\beta\}$ are transcendence bases of the field of rational functions of
$\mathcal{V}$. It follows that $\mathcal{I}(V) \cap k[\gamma]=0$, that
is, the set of possible tuples of residual coefficients of the points $\widetilde{q}_i$
contains a dense Zariski open set.
\end{proof}

\begin{exmp}
Consider the case of two conics $C_1=`` (-11) +2x +2y +2xy +0x^2 +0y^2"$,
$C_2=`` 0 +8x +14y +20xy +12x^2 +14y^2"$, their stable intersection is the set
of points $\{(2, -6), (-4,2), (-13,-14), (-6,-6)\}$. These four points are in
general position with respect to $C_1$ and $C_2$ so, for any generic lifts of
$C_1$, $C_2$, the residual coefficients of their intersection points are
generic. However, consider now the case of two conics $C_1= ``0 +(-10)x +(-10)y
+(-10)xy +0x^2 +0y^2"$ and $C_2=``0 +(-10)x +(-10)y +(-10)xy +1x^2 +2y^2"$. They
have only one intersection point of multiplicity $4$, $\{(0,\frac{1}{2})\}$. Taking the point three or four times yields to a set which is not in general position in none of the curves. Hence, the maximal number of intersection points that are in general position in both curves is 2. So, the drawback of this theorem is that the number $n$ of points in general position in both curves is not uniform with respect to the supports. The following is a uniform result that holds for every pair of curves with prescribed support.
\end{exmp}

\begin{thm}\label{un_punto_interseccion_es_generico}
Suppose given two tropical curves $C_1$, $C_2$ with support $I_1$ and $I_2$ respectively. Let $\widetilde{C}_1$, $\widetilde{C}_2$ be two lifts of the curves whose principal coefficients are generic and let $q$ be one stable intersection point. Then, the principal residual coefficients of $\widetilde{q}$ are generic. That is, if we impose polynomial conditions $F\neq 0$ to the coefficients of $\widetilde{C}_i$ then the possible residual coefficients of the point $\widetilde{q}$ contains a dense constructible set of $k^2$.
\end{thm}
\begin{proof}
One point $q$ is always in general position with respect to any curve, so we are
in the hypothesis of
Theorem~\ref{teorema_de_puntos_interseccion_en_posicion_general}
\end{proof}

\section{Lift of a Construction}\label{section_lift}

Let $\mathfrak{C}$ be a geometric construction of graph $G$. This Section
deals with the problem of lifting a tropical instance of $G$ obtained by the
construction to an algebraic instance. Let $H_0$ be the set of input elements of
$\mathfrak{C}$ and $h$ a tropical realization of $H_0$. The steps of the
construction define a tropical realization $p$ of $G$. On the other hand, let
$\widetilde{h}=T^{-1}(h)$ be any algebraic realization of $H_0$ that projects
onto $h$ (recall that this lift is not unique). Then, there are two potential
problems. First, it is possible that $\mathfrak{C}$ is not well defined in
$\widetilde{h}$. Second, if the construction is well defined and $\widetilde{p}$
is the algebraic realization of $G$ obtained from $\widetilde{h}$, it is
possible that $T(\widetilde{p})\neq p$. We study
conditions for the lift $T^{-1}(h)$ such that the following Diagram commutes:
\begin{equation}\label{diagrama_tropical}
\begin{matrix} (\mathbb{K}^*)^2&&\mathbb{T}^2\\ \text{
Input}\ \widetilde{h}&\stackrel{T^{-1}}{\longleftarrow}&\text{Input}\ h\\
\mathfrak{C}\downarrow&&\downarrow\mathfrak{C}\\ \text{
Output}\ \widetilde{p}&\stackrel{T}{\longrightarrow}&\text{
Output}\ p
\end{matrix}
\end{equation}

Given an instance of a geometric construction, we define
sufficient residual conditions on the lifts $\widetilde{h}$ of the input $h$ for
the compatibility $T(\widetilde{p})=p$. In order to do this, let $\{C_1, \ldots,
C_n, q_1, \ldots, q_m\}$ be the input elements of a geometric construction
$\mathfrak{C}$, curve $C_i$ of support $I_i$, point $q_j \in (\mathbb{T}^*)^2$.
Take $N=2m +\sum_{i=1}^n (\delta(I_i)-1)$ and let $\{\widetilde{f}_1, \ldots,
\widetilde{f}_n, \widetilde{q}_1, \ldots, \widetilde{q}_m\}$ be a set of lifts
of a concrete tropical instance of the input, $f_i = \sum_{(k,l) \in I_i}
\widetilde{a}_{(k,l)}^i x^k y^l$, $\widetilde{q}_j = (\widetilde{q}_j^1,
\widetilde{q}_j^2)$. We are going to compute a constructible set $\mathfrak{S}
\subseteq (k^*)^N$, not always empty, that encodes the residual conditions for
the compatibility of the algebraic and tropical construction. We are going to
define two auxiliary sets $T$ and $V$ first. The set $T$ is defined adding the
residual restrictions obtained by Theorems~\ref{condiciones_de_pseudocramer} and
\ref{interestable} that ensure that each step of the construction is compatible
with tropicalization. Let
\[f_i=``\sum_{(k,l) \in I_i} a^i_{(k,l)} x^k y^l", 1\leq i\leq n,\]
\[q_j= (q_j^1, q_j^2), 1\leq j\leq m\]
be the tropical input elements. Take a generic lift of the input
\[\widetilde{f}_i'=\sum_{(k,l) \in I_i} \widetilde{a}^i_{(k, l)} x^k
y^l, 1\leq i\leq n,\]
\[\widetilde{q}_j'=(\widetilde{q}_j^1, \widetilde{q}_j^2 ), 1\leq j\leq
m\]
and $V_0= \{ \alpha^i_{(k, l)}, \gamma^r_j\}$ is a set of indeterminates where
$Pc(\widetilde{a}^i_{(k,l)})=\alpha^i_{(k,l)}$,
$Pc(\widetilde{q}^i_j)=\gamma_j^i$. These indeterminates will describe
$\mathfrak{S}$. Perform the construction with this data as follows.

Start defining the constructible set $T=(k^*)^N= \{x \in k^N | \alpha^i_{(k, l)}
\neq 0, \gamma^r_j\neq 0, 1\leq i\leq n, 1\leq j\leq m\}$ and $V=V_0$. We are
going to redefine $T$ and $V$ inductively at each step of the construction.
Suppose that we have defined $V$ and the constructible set $T\subseteq (k^*)^V$
for the construction up to a construction step. We redefine $T$ after the step
as follows: For the case of the computation of the curve $C$ of support $I$
passing through $\delta(I)-1$ points, we have to solve a system of linear
equations. The coefficients of $\widetilde{C}$ are rational functions of the
variables $V$. Theorem~\ref{condiciones_de_pseudocramer} provides sufficient
conditions in the variables $V$ for the system being compatible with
tropicalization. These conditions are of the form $\Delta_{A^i}(Pc(\widetilde{A}^i))\neq 0$ where
$A$ is the tropical matrix of the system of linear equations. We add to $V$
$(\delta(I)-1)$ new variables $s_1, \ldots, s_{\delta-1}$ and we consider
$T\subseteq (k^*)^{K+ \delta -1}$. We add the conditions
$\Delta_{A^i}(Pc(\widetilde{A}^i)) \neq 0$ to the definition of $T$ and the
equations $\Delta_{A^i} (Pc( \widetilde{A}^i)) - s_i \Delta_{ A^{i_0}} (Pc(
\widetilde{A}^{ i_0 })) =0$, where $i_0$ is a dehomogenization variable of $C$.
We follow the construction with $\widetilde{C}$ among our available objects.

Suppose now that our construction step consists in the intersection of two
curves $\widetilde{f}$, $\widetilde{g}$ of support $I_f$, $I_g$ respectively.
Its stable intersection can be determined using the technique of resultants.
That is, consider first the
resultant polynomials $\widetilde{R}_x(x)= \text{Res}( \widetilde{f}, \widetilde{g},y)$,
$\widetilde{R}_y(y)= \text{Res}( \widetilde{f}, \widetilde{g},x)$. Let $a$ be a natural
number such that $x-ay$ is injective in the finite set $\mathcal{T}(f) \cap
\mathcal{T}(g)\cap \mathcal{T}(R(x)) \cap \mathcal{T}(R(y))$. Let
$\widetilde{R}_z(z) = \text{Res} (\widetilde{f} (zy^a,y), \widetilde{g}
(zy^a,y),y)$. If $t_r$ are the variables of $V$ corresponding with the principal
coefficients of $\widetilde{f}$, $\widetilde{g}$, Theorem~\ref{interestable}
provides sufficient conditions of the form $\widetilde{u}(t_r)\neq 0$ that
ensures that the algebraic and tropical intersection are compatible. We add
these polynomials $\widetilde{u}(t_r)\neq 0$ to the definition of $T$. In the
tropical context, there are $M=\mathcal{M} (\Delta_f, \Delta_g)$ stable
intersection points $b_j=(b_j^1,b_j^2)$. We add $2M$ new variables $s_j^1$,
$s_j^2$, $1\leq j\leq M$ to $V$. Consider $T$ contained in $(k^*)^{K+2M}$. For
each tropical point $b_j$, let $s_{j_1}, \ldots, s_{j_n}$ be the algebraic
points projecting into $b_j$. We take the following equations:
\[(\widetilde{R}_x )_{ b_j^1} = \prod_{ r=1}^n (x -s_{j_r}^1), \qquad
(\widetilde{R}_y)_{ b_j^2}
= \prod_{r=1}^n (y -s_{j_r}^2),\]
\[(\widetilde{R}_z)_{``b_j^1(b_j^2)^{-a}"}=\prod_{r=1}^n
(z-s_{j_r}^1(s_{j_r}^2)^{-a}).\]
In this way, the coefficients of $(\widetilde{R}_x)_{b_j^1},
(\widetilde{R}_y)_{b_j^2}$ and $(\widetilde{R}_z)_{``b_j^1(b_j^2)^{-a}"}$ are
identified with symmetric functions in $s_{j_r}^1$, $s_{j_r}^2$ and
$s_{j_r}^1(s_{j_r}^2)^{-a}$ respectively. We add these identifications to the
definition of $T$. In this way, we ensure that there is a bijection between the
roots of the resultants and the variables $s_j$. We also add the residual
conditions of the curves over the intersection points
$\widetilde{f}_{b_j}(s_j^1, s_j^2)=0$, $\widetilde{g}_{b_j}(s_j^1, s_j^2)=0$,
and the conditions of the points being in the torus $s^1_js^2_j\neq 0$. We
continue the construction with the points $(s_i^1 t^{-b^1_i}, s_i^2
t^{-b_i^2})$. Notice that we are only defining the principal terms of the
elements, because this is all the information needed for the Theorem. After the
whole construction, we have defined a constructible set $T$ that characterizes
the possible principal term of every element in the construction. Finally,
$\mathfrak{S}$ is defined as the projection of the set defined by $T$ into the
space of variables $V_0.$

\begin{defn}\label{conjunto_S}
The set $\mathfrak{S}$ previously defined is called \emph{the set of valid
principal coefficients of the input elements}.
\end{defn}

\begin{thm}\label{versionconstruct}
Let $\{C_1, \ldots, C_n, q_1, \ldots, q_m\}$ be the input elements of a
geometric construction $\mathfrak{C}$, curve $C_i$ of support $I_i$, point $q_j
\in (\mathbb{T}^*)^2$. Take $N=2m +\sum_{i=1}^n (\delta(I_i)-1)$ and let
$\{\widetilde{f}_1, \ldots, \widetilde{f}_n, \widetilde{q}_1, \ldots,
\widetilde{q}_m\}$ be a set of lifts of a concrete tropical instance of the
input, $\widetilde{f}_i = \sum_{(k,l) \in I_i} \widetilde{a}_{(k,l)}^i x^k y^l$, $\widetilde{q}_j = (\widetilde{q}_j^1, \widetilde{q}_j^2)$, $Pt(\widetilde{a}^i_{(k,l)})=\alpha^i_{(k,l)}t^{-a^i_{k,l}}$, $Pt(\widetilde{q}_j^i)=\gamma_j^it^{-q_j^i}$. Let $\mathfrak{S} \subseteq (k^*)^N$ be the set of valid principal coefficients of the input. Then, if the vector
\[\big(\alpha^1_{ (k, l) }, \ldots, \alpha^n_{ (k, l) }, \gamma^1_1, \ldots,
\gamma_m^2\big) \in (k^*)^N\]
of principal coefficients lies in $\mathfrak{S}$, the algebraic construction is
well defined and the result projects onto the tropical construction.
\end{thm}
\begin{proof}
Suppose that the vector $(\alpha^1_{(k,l)}, \ldots, \alpha^n_{(k,l)},
\gamma^1_1, \ldots, \gamma_m^2)$ belongs to $\mathfrak{S}$. We are going to
construct suitable algebraic data. Perform the steps of the construction. For
the curve passing through a number of points, the set $\mathfrak{S}$ imposes
that there is only one solution of the linear system we have to solve and that
this solution projects correctly. For the case of the intersection of two
curves, the resultants $\widetilde{R}_x$, $\widetilde{R}_y$, $\widetilde{R}_z$ are compatible with projection. So, the curves intersects in finitely many points in the torus and these points projects correctly onto the tropical points. So this step is also compatible with the tropicalization.
\end{proof}

In this theorem, it is not claimed that there is always a possible lift, as
Theorem~\ref{versionincidente_aciclico} does. It is possible that the set
$\mathfrak{S}$ is empty. In this case, the theorem do not yield to any
conclusion. In Subsection~\ref{certificates} we will discuss what
can be said if $\mathfrak{S}$ is empty.

Now, we search sufficient conditions for a
construction $\mathfrak{C}$ that assert that the set $\mathfrak{S}$ is non empty
for every realization $h$ of the input. For example, let $\mathfrak{C}$ be a
depth 1 construction. There are only two kind of elements, input elements and
depth 1 elements. If the realization $\widetilde{h}$ of the input elements is
generic, by Theorems~\ref{condiciones_de_pseudocramer} and \ref{interestable},
every depth 1 element is well defined and projects correctly. Thus, every depth
1 tropical construction can be lifted to the algebraic plane. Furthermore, if the vector
of coefficients of the depth 1 elements is generic, we would be able to
construct some other depth 2 elements from them. By
Theorems~\ref{la_curva_por_puntos_genericos_es_generica} and
\ref{un_punto_interseccion_es_generico}, we already know that every single depth
1 element is generic. However, it may happen that there are algebraic relations
among the set of depth 1 elements that do not allow to apply induction in
further steps. So, in order to use an induction scheme over the construction, we
need to ensure that in future steps of the construction we will only use
elements that are generic. Next Definition describes constructions such that
this genericity of the elements always holds, whatever the input elements are.

\begin{defn}\label{def_geo_admissible}
Let $\mathfrak{C}$ be a geometric construction. Let $G$ be the incidence graph
with the orientation induced by the construction. The construction
$\mathfrak{C}$ is \emph{admissible} if, for every two nodes $A$, $B$ of $G$,
there is at most one oriented path from $A$ to $B$. In the case where the
construction is not admissible, let $A$, $B$ two elements such that there is at
least two paths from $A$ to $B$. This is denoted by $A \rightrightarrows B$.
\end{defn}

The main Theorem of the Section proves that if $\mathfrak{C}$ is an
admissible geometric construction, then every tropical realization of
$\mathfrak{C}$ can be lifted to a compatible algebraic realization.

\begin{thm}\label{admisible_es_levantable}
Let $\mathfrak{C}$ be an admissible geometric construction. Then, for every
tropical instance of the construction, the set $\mathfrak{S}$ defined in
Theorem~\ref{versionconstruct} is nonempty and dense in $(k^*)^N$. Moreover, for
every element $X$ of the construction, its possible values, as the input
elements range over $\mathfrak{S}$, contains a dense open subset of its support
space. In particular, every tropical instance of the construction $\mathfrak{C}$
can be lifted to the algebraic plane $(\mathbb{K}^*)^2$.
\end{thm}

\begin{proof}
We prove the Theorem by induction in the depth of the construction. If the
construction is of depth 0, then there is nothing to prove, because the set of
steps is empty and $\mathfrak{S}=(k^*)^N$ which is dense and the values of each
element are dense in their respective space of configurations. Suppose the
Theorem proved for admissible constructions of depth smaller or equal to $i$.
Let $\mathfrak{C}$ be any admissible construction of depth $i+1$. For each
element $X$ of depth $i+1$, let $Y_1, \ldots, Y_n$ be the direct predecessors of
$X$. By induction hypothesis, the set of possible values of $Y_i$ contains a
dense open set in its space of configurations. As the construction is
admissible, the set of predecessors of $Y_i$ is disjoint from the set of
predecessors of $Y_j$, if $i\neq j$. Because if both elements had a common
predecessor $A$, there would be a double path $A \rightrightarrows X$, contrary
to the hypothesis. Hence, the coefficients $Y_1,\ldots,Y_{n}$ are completely
independent and the possible tuples ($Y_1$, \ldots, $Y_n$) are just the
concatenation of possible values of coefficients of each element $Y_i$. By the
results in the Theorems~\ref{la_curva_por_puntos_genericos_es_generica} and
\ref{un_punto_interseccion_es_generico}, as the elements $Y_j$ are generic, so
is $X$. That is, the possible values of $X$ contains a dense open set of its
support space. The conditions imposed by the definition of $X$ to the auxiliary
set $T$ in Theorem~\ref{versionconstruct} are a set of inequalities in the
tuples $(Y_1,\ldots,Y_n)$ that are verified on an open set. Likewise, the
restrictions in the elements $Y_j$ impose other restrictions to their
predecessors. Again, this restrictions are verified in an open set, we are
explaining this with more detail:

If $Y_j$ is constructed from elements $Z_{ji}$, there is a set of restrictions
$f_{s}(Z_{ji})\neq 0, s\in S$ that ensure that $Y_j$ is well defined and it is
compatible with tropicalization. Let $g_l(Y_1,\ldots,Y_n)\neq 0$, $l\in L$ be
the polynomials imposed by $X$ to be well defined and compatible with
tropicalization. In addition to this, if $Y_j=(Y_j^1,\ldots,Y_j^{n_j}),$ each
variable $Y_{j}^r$ is algebraic over the field $p(k)(Z_{ji})$, where $p(k)$ is
the prime field of $k$. If we multiply each polynomial $g_l(Y_1,\ldots,Y_n)$ by
its conjugates in the normal closure of $p(k)(Z_{ji})\subseteq
p(k)(Z_{ji},Y_{i})$, we obtain some polynomials $G_l(Z_{j1},\ldots,Z_{jn})$. If
neither $G_l(Z_{ij})$ nor $f_{s}(Z_{ij})$ are zero, then the elements $Y_i$ and
$X$ are well defined and are compatible with projection. These polynomials
define possible valid principal coefficients for the subconstruction
$Z_{ji}\rightarrow Y_i\rightarrow X$. Applying this method recursively, we
obtain a set of conditions in the input elements. Let $\mathfrak{S}_i$ the set
of good input elements for every subconstruction of $\mathfrak{C}$ consisting on
the elements of depth up to $i$. By induction hypothesis, $\mathfrak{S}_i$ is
non empty and contains an open Zariski set. Intersecting this set with the open
sets induced by each element $X$ of depth $i+1$ to be compatible with
tropicalization, we obtain that the required set $\mathfrak{S}_{i+1}$ contains a
dense open Zariski set.
\end{proof}

Tropical geometric constructions are a useful tool when dealing with non-trivial
incidence relations between varieties. It agrees naturally with the stable
intersection of the curves taken in consideration. Moreover, it permits to
arrange the computations focusing on the smaller set of input objects. Now, we quantify how well a realization of a construction behave with respect to tropicalization. In order to determine the potentially good situations, we focus on the following concepts:
\begin{itemize}
\item An abstract geometric construction. That is, we do not specify the coordinates of the points, neither the concrete curves, only their support and the steps of the construction. Moreover, we ask it to be well defined in both fields $\mathbb{K}$ and $k$.
\item The specialization of the input elements of the abstract construction to
concrete tropical elements.
\item A concrete algebraic lift of the given set of input elements.
\end{itemize}

These concepts are manipulated by adding quantifiers relating them in order to
obtain a statement like:\\
``\emph{$K_1$ tropical construction $K_2$ specialization of the input data $K_3$
lift of these input data, diagram~\ref{diagrama_tropical} commutes}".\\
Where $K_1$, $K_2$, $K_3\in\{\forall,\ \exists\}$. We arrive naturally to the
following problems:

\begin{prob}\label{questions}\mbox{}
\begin{enumerate}
\item For all constructions, for all input tropical data and for all lifts of
these tropical data, diagram~\ref{diagrama_tropical} commutes.
\item For all constructions and for all input tropical data there exists a lift
of these tropical data such that diagram~\ref{diagrama_tropical} commutes.
\item For all construction, there is a choice of the input tropical data such
that for all lift of these tropical data, diagram~\ref{diagrama_tropical}
commutes.
\item There exists a construction such that for all input tropical data and for
all lifts of these tropical data, diagram~\ref{diagrama_tropical} commutes.
\item For all constructions, there is a choice of input tropical data and there
is a lift of these tropical data such that diagram~\ref{diagrama_tropical}
commutes.
\item There exists a construction such that for all input tropical data there
is a lift of these tropical data such that diagram~\ref{diagrama_tropical}
commutes.
\item There exists a construction and there is suitable input tropical data
such that for all lifts of these tropical data, diagram~\ref{diagrama_tropical}
commutes.
\item There exists a construction, particular input tropical data and a suitable
lift of these tropical data such that diagram~\ref{diagrama_tropical} commutes.
\end{enumerate}
\end{prob}

Clearly, these relations are not independent, ranking (non linearly) from item
1, which is the strongest, to item 8, the weakest one. Checking this problems
gives an overview of the typical problems we find when dealing with incidence
conditions in Tropical Geometry. The only statements that hold are items 5, 6, 7
and 8. For the sake of brevity, we will consider mostly the case where our
curves are lines on the plane.

\begin{prop}
The only items of problem~\ref{questions} that hold are 5, 6, 7 and 8.
\end{prop}
\begin{proof}\mbox{}

$\bullet$ Take two tropical lines in the plane that intersects in only one
point. Then, for all lifts of this two lines, the intersection point always
tropicalizes to the tropical intersection. So statement \ref{questions}.7 holds
and, from this, we derive that \ref{questions}.8 also does.

$\bullet$ Choose two curves that intersect in an infinite number of points. In
Theorem~\ref{versionincidente_aciclico}, we are given a way to compute lifts
that
intersects in non stable points. So the property of agreement with
tropicalization is not universal for the non transversal cases. This simple
example shows that statement \ref{questions}.1 does not hold. Using duality, we
observe also that the concept of stable curve through a set of points does not
work for every input data and every lift (ie. there will always be exceptional
cases). Thus, since every tropical geometric construction consists of a sequence
of these two steps (computing the stable curve through a set of points, or
computing the stable intersection of two curves), we deduce that statement
\ref{questions}.4 neither holds. In particular, if we are able to find a
construction such that for all input data we arrive to these exceptional cases,
we
will find a counterexample to question \ref{questions}.3. An example of such a
construction is as follows:\\
\begin{tabular}{ll}
Input: & points $a$, $b$, $c$, $d$, $e$.\\
Depth 1: & lines $l_1:=\overline{ab}$, $l_2:=\overline{ac}$,
$l_3:=\overline{ad}$, $l_4:=\overline{ae}$.\\
Depth 2: & points $p_{12}=l_1\cap l_2$, $p_{13}=l_1\cap l_3$, $p_{14}=l_1\cap
l_4$,
$p_{23}=l_2\cap l_3$,\\
 & $p_{24}=l_2\cap l_4$, $p_{34}=l_3\cap l_4$.
\end{tabular}

First, we compute four tropical lines through one fixed point $a$. If point $a$
is exactly the vertex of one of the lines, then two of the input points are the
same and there is an infinite number of lines passing through these two points.
On the other hand, if $a$ is never the center of the lines, it must be in one of
the three rays. There are only three possibilities for the rays, the directions
$(-1,0)$ $(0,-1)$ and $(1,1)$. As there are four lines involved, two of the
branches must have the same directions, so these two lines intersect in an
infinite number of points and we are done.

$\bullet$ To go further in the analysis, it is necessary to have more tools that
takes care of more complicated constructions.
Theorem~\ref{admisible_es_levantable} establishes that for an admissible
construction and for all realization of the input elements, there always exists
a lift of these elements such that all the steps of both constructions are
coherent with the tropicalization. In particular, we have the validity of
question \ref{questions}.6 for every admissible construction.

$\bullet$ Also, a counterexample to \ref{questions}.2 is the following. Take three points $a$, $b$, $c$. Construct the lines $l_1=\overline{ab}$, $l_2=\overline{ac}$ and the point $p=l_1\cap l_2$. If we perform this construction in the projective plane with three points not in the same line, we will always find that $p=a$. But in the tropical case, taking $a=(0,0)$, $b=(-2,1)$, $c=(-1,3)$, we arrive to $p=(0,1)\neq a$. This simple example shows a concrete construction and input data such that for all lifts of the input elements, diagram \ref{diagrama_tropical} does not commute. Note that in this case there are double paths in the construction graph. If we follow the method exposed in Theorem~\ref{versionconstruct}, then, for all lifts, we arrive that the constructible set $\mathfrak{S}$ is contained in $0\neq 0$. That is, the set of valid principal coefficients is empty.

$\bullet$ Finally, let us prove \ref{questions}.5. This case of course cannot be
restricted to the linear case. Suppose given a geometric construction, we choose
as input data the most degenerate case possible: if we have a point, we choose
the point to be $p_0:=(0,0)$ and if we have a curve with prescribed support, we
take all its coefficients equal to zero. As a set, it consists in some rays
emerging from the origin $(0,0)$ in perpendicular directions to the edges of the
Newton polygon of the curve. The stable intersection of any two such curves is
always the isolated point $p_0$ with the convenient multiplicity. The stable
curve with prescribed support taking all elements equal to the origin is the
one with all coefficient equal to zero. It only rests to check that there is a
lift compatible with this tropical construction. As the construction is well
defined, it is realizable for the generic input in $(k^*)^2$. This construction
can be embedded in $(\mathbb{K}^*)^2$ with all the elements of order 0.
\end{proof}

As an application of the construction method and
Theorem~\ref{admisible_es_levantable}, we are able to extend
Theorem~\ref{versionincidente_aciclico} to a wider set of incidence
configurations.

\begin{thm}
Let $G$ be an incidence structure, suppose that we have a tropical realization
$p$ of $G$ such that, for every curve $C$, the set of points incident to $C$ are
in generic position with respect to $C$. Then, the tropical realization can be
lifted to an algebraic realization.
\end{thm}
\begin{proof}
For each curve $C$ of support $I$, let $q_1,\ldots, q_n$ be the set of points
incident to $C$. By definition of points in general position, we can extend
this set to a set of points $q_1, \ldots, q_{\delta(I)-1}$ such that $C$ is the
stable curve through these points. Add to the configuration $G$ these additional
points for every curve $C$. We obtain in this way an incidence configuration
$G_1$ that contains $G$ as a substructure and such that every curve $C$ of
support $I$ is exactly the stable curve passing through the points $q_1, \ldots,
q_{\delta(I)-1}$. Hence, by Proposition~\ref{subgrafodeconstruccion}, $G_1$ is
the graph of a geometric construction $\mathfrak{C}$. The input elements are the
set of points $q_i$ and every curve is the stable curve through $\{q_1, \ldots,
q_{\delta(I)-1}\}$. This construction is admissible, because it is of depth 1.
By Theorem~\ref{admisible_es_levantable}, every tropical instance $p_1$ of
$\mathfrak{C}$ can be lifted to an algebraic instance $\widetilde{p}_1$ of
$\mathfrak{C}$. In particular, the instance $p$ of $G$ we started from can be
lifted to the algebraic plane.
\end{proof}

This Theorem shows how the notion of points in general position helps to the
problem of lifting an incidence configuration. Our next goal is to apply this
notion to more complex configurations coming from geometric constructions. The
key idea for this application is that points in general position with respect
to a curve $C$ behave like generic points for the purposes of
Theorem~\ref{admisible_es_levantable}.

\begin{thm}\label{Construccion_casi_admisible}
Suppose that we are given a non admissible geometric construction $\mathfrak{C}$
but such that the only obstacle to be an admissible construction is that we
have two curves $C_1$, $C_2$ with intersection $Q=\{q_1, \ldots, q_n\}$ such
that $Q$ is used twice to define some successor element $x$. That is, every
double path $A\rightrightarrows B$ in $\mathfrak{C}$ can be restricted to a
double path from both curves passing through $Q$, \[C_1 \rightrightarrows
Q\rightrightarrows B\ \textrm{and\,}\ C_2 \rightrightarrows Q\rightrightarrows
B.\] Suppose we have an instance $p$ of this construction. If, for every element
$x$ which is the end of a double path, the set $Q_x = \{q_i\in Q\ |\ \exists q_i
\rightarrow x \}$ is in general position in $C_1$ and $C_2$, then the tropical
instance can be lifted to an algebraic realization $\widetilde{p}$ of the
construction. More concretely, the set $\mathfrak{S}$ of
Theorem~\ref{versionconstruct} associated to $p$ contains an open dense subset
of $(k^*)^N$.
\end{thm}
\begin{proof}
First, we are proving that, for any single node $x$ of $\mathfrak{C}$, its
construction can be lifted. Let $x$ be a node of $\mathfrak{C}$. Let
$\mathfrak{C}_x$ be the minimal subconstruction of $\mathfrak{C}$ such that it
contains every input element of $\mathfrak{C}$ and the element $x$. This minimal
subconstruction can be defined as follows. First, we consider as nodes of
$\mathfrak{C}_x$ the input elements of $\mathfrak{C}$, the node $x$ and every
predecessor of $x$. The incidence conditions will be those induced by
$\mathfrak{C}$. Second, we complete it with the necessary nodes of
$\mathfrak{C}$ as in the proof of Proposition~\ref{subgrafodeconstruccion}.
Actually, the only nodes we have to add are the intersection points of two
curves $y_1$, $y_2$ that have to be intersected (necessarily, these curves will
be predecessors of $x$). Let $\mathfrak{S}_x$ be the set of valid input elements
of the construction $\mathfrak{C}_x$. By the construction of $\mathfrak{S}$,
\[\mathfrak{S}=\bigcap_{x\in \mathfrak{p}\cup \mathfrak{B}} \mathfrak{S}_x\] So,
if every $\mathfrak{S}_x$ contains a non empty open Zariski set of $(k^*)^N$,
the same occurs for $\mathfrak{S}$.

If $\mathfrak{C}_x$ is admissible, then $\mathfrak{S}_x$ contains a non empty
Zariski set by Theorem~\ref{admisible_es_levantable}. If $\mathfrak{C}_x$ is not
admissible, the set $Q_x$ contains at least two elements. Moreover, for every
node $y$ in $\mathfrak{C}_x$ it happens that $Q_y \subseteq Q_x$.

Consider now the minimal subconstruction $\mathfrak{C}_x^1$ containing every
input element and the set $Q_x$. This construction is admissible, so
$\mathfrak{S}_x^1$ is dense. On the other hand, the possible principal
coefficients of the set $Q_x$ form a dense set of its space of configurations by
Theorem~\ref{teorema_de_puntos_interseccion_en_posicion_general}. Let
$\mathfrak{C}_x^2$ the subconstruction obtained from $\mathfrak{C}_x$ by
deleting every predecessor of the points in $Q_x$ and the intersection of $C_1$
and $C_2$ not in $Q_x$. This construction is also admissible, because the curves
$C_1$, $C_2$ have been deleted among other objects. Hence $\mathfrak{S}_x^2$ is
also dense. The projections of the set $\mathfrak{S}_x^1$ and $\mathfrak{S}_x^2$
into the support space of $Q_x$ contains an open dense subset, their
intersection also contains a non empty dense subset. This means that there are
values of the principal coefficients of $Q_x$ that are generic and
compatible either with $\mathfrak{C}_x^1$ and $\mathfrak{C}_x^2$. It follows
that for a residually generic lift of the input elements of $\mathfrak{C}_x$,
every step will be well defined and compatible with tropicalization. Thus,
$\mathfrak{S}_x$ is contains a dense subset of $(k^*)^N$.
\end{proof}

In contrast to Theorem~\ref{admisible_es_levantable}, this Theorem does
not work for every tropical realization of a particular construction $\mathfrak{C}$,
because it is stated in terms of the realization. It needs some additional
hypothesis in the construction (some points are in general position) that depend
on the concrete realization. It still has its applications, such as
Theorem~\ref{weak_pascal}.

\subsection{Impossibility for the Existence of a Lift}\label{certificates}
Now we work with non admissible constructions, suppose that we
have a non admissible geometric construction $\mathfrak{C}$ and a tropical
instance of it such that the constructible set $\mathfrak{S}$ is empty. Then, we
would still like to know if it is possible to lift the construction. The only
result that affirms that it is impossible to have a lift is
Proposition~\ref{certificado_curva_por_puntos}. We can provide a similar notion
for the stable intersection of curves. Theorem~\ref{interestable} provides
compatibility restrictions in the residual coefficients of $\widetilde{f}$ and $\widetilde{g}$
in terms of resultants $R(x)$. Next proposition states some certificates of
compatibility and incompatibility between the algebraic and the tropical
resultant.

\begin{prop}\label{certificado_interseccion}
Let $f,g$ be two tropical curves, let $\{\gamma_1, \ldots, \gamma_{r}\}$ be the
residual conditions for the compatibility of the algebraic and tropical
resultant $R(x)$ provided by Theorem~\ref{interestable}. These are the residual
coefficients $\gamma_i\neq 0$ corresponding to the indices $i$ that are vertices of the subdivision induced in the Newton polygon of $R(x)$. With these conditions:
\begin{itemize}
\item If every polynomial $\gamma_i$ is a monomial, then, the algebraic
resultant is always compatible with tropicalization
$T(\widetilde{R}(x))=\mathcal{T}(R(x))$.
\item If one polynomial $\gamma_i$ is a monomial, then the algebraic
resultant $\widetilde{R}(x)$ is compatible with tropicalization if and only if
the rest of the polynomials $\gamma_j$ are non zero.
\item If every polynomial $\gamma_i$ is zero, we cannot derive any information
about the compatibility.
\end{itemize}
\end{prop}
\begin{proof}
Let $\widetilde{R}(x)=\sum_{i=0}^r \widetilde{h}_ix^i$, $R(x)=\sum_{i=0}^r,
h_ix^i$ be the algebraic and tropical resultant. If $\gamma_i\neq 0$ then the
principal term of $\widetilde{h}_i$ is exactly $\gamma_it^{-h_i}$. The
conditions searched for the compatibility of the resultants is that the elements
$\gamma_i$ associated to an index $i$ such that it is a vertex of the subdivision induced in the Newton polytope of $R(x)$ do not vanish. If one $\gamma_i$
is a monomial, then it will never evaluate to zero. So the Newton diagram will
not change if and only if the rest of the $\gamma_j$ do not evaluate to zero.
Hence we have the first two items. On the other hand, if every $\gamma_i$
evaluates to zero, we cannot know how the Newton diagram of $\widetilde{R}(x)$
is with respect to the Newton diagram of $R(x)$, it may change or not.
\end{proof}

\begin{defn}
Let $\mathfrak{C}$ be a construction and $p$ a tropical realization of it. Let
$x$ be a node of $\mathfrak{C}$. We say that $x$ is a \emph{fixed} element of
$\mathfrak{C}$ if:
\begin{itemize}
\item $x$ is an input element of $\mathfrak{C}$.
\item $x$ is the curve of support $I$ passing through $\{y_1,\ldots,
y_{\delta(I)-1}\}$ and at least one of the tropical minors of the linear system
defining $x$ is regular (See Proposition~\ref{certificado_curva_por_puntos}).
\item $x$ is an intersection point of $y_1$ and $y_2$ and, if $C_1$, $C_2$
are the tropical realization of curves $y_1$, $y_2$, then, at least one the
residual conditions $\gamma_{i_1}(x)$, $\gamma_{i_2}(y)$ and
$\gamma_{i_{3}}(xy^{-a})$ of each resultant $R(x)$, $R(y)$, $R(xy^{-a})$ defined
in Theorem~\ref{interestable} is a monomial.
\end{itemize}
\end{defn}

Let $\mathfrak{C}$ be a geometric construction and $p$ a tropical realization of
$\mathfrak{C}$. Suppose that the set $\mathfrak{S}$ associated to the tropical
realization is empty. Then, during the definition of the auxiliary set $T$ in
Theorem~\ref{versionconstruct}, there will be a step such that $T$ was not
empty before the step, but the restrictions added in this step forces $T$ to be
empty. This step consists in defining an element $x$. Let $h_1, \ldots, h_r$ be
the residual polynomials codifying the compatibility of this algebraic step with
tropicalization defined using Theorem~\ref{condiciones_de_pseudocramer} and
\ref{interestable}. Suppose that at least one of the polynomials $h_i$ does not
evaluate to zero. Then:
\begin{itemize}
\item If every predecessor of $x$ is fixed, by
Propositions~\ref{certificado_curva_por_puntos} and
\ref{certificado_interseccion}, there cannot be any lift of the tropical
realization of $\mathfrak{C}$. Because for every lift of the input elements,
either one of the predecessors of $x$ does not tropicalize correctly or, if every lift of the predecessors of $x$ tropicalize correctly, then the element $x$ either is not well defined, or it will never tropicalizes correctly.
\item If at least one predecessor of $x$ is not fixed, then, there might be a lift of the tropical realization of $\mathfrak{C}$ or not. But at least, there cannot be any lift with residually generic input elements. There must be some algebraic relations among the residual coefficients of the algebraic input elements of $\mathfrak{C}$.
\end{itemize}

On the other hand, if every residual polynomial $h_i$ evaluates to zero. We
cannot conclude anything, there might be a lift of the realization or not. And
this lift may work for the generic input or not. In this case the residual
coefficient approach is not enough to answer the question.

For most geometric constructions the remarks above are enough. That is, if for
one tropical realization its associated set $\mathfrak{S}$ is empty, then either
we can deduce that for the generic lift of the input elements the algebraic
construction will not project correctly. Or even that there will be no lift at
all. In fact, for every geometric construction that we have faced during the
development of this theory, every instance of every construction fell in these
two cases. It is difficult to find a construction and an instance of the
construction such that the construction method and the set $\mathfrak{S}$ does
not provide any information. The following example is basically the only one with this behavior that we are aware of.

\begin{exmp}
In this example, for convenience with the geometric language, we will think that
the algebraic torus $(\mathbb{K}^*)^2$ is contained in the affine plane and this one
contained in the projective plane. With this in mind, we can talk about concepts
such at horizontal line (curve of support $\{(0,0),(0,1)\}$) vertical line (curve of support $\{(0,0),(1,0)\}$) or the line at the infinity. This is intended only to simplify notations and use a more natural language, but it does not interfere with the result itself.

First, we need a specific construction. Given a point $a$ and a line $l$. We
look for a geometric construction such that, in the algebraic plane, it defines
the parallel of $l$ passing through $a$. The difficulty is to define it with the
restricted allowed steps of Definition~\ref{construccion}.

\noindent
\begin{tabular}{ll}
\multicolumn{2}{l}{$l'$=Parallel($a$, $l$, $q$):}\\
Input: &points $a$, $q$, line $l$.\\
Depth 1: & vertical line $v_1$ passing through $a$.\\
 &vertical line curve $v_2$ passing through $q$.\\
 &horizontal line $h_1$ passing through $q$.\\
 &line $r_1$ passing through $\{a,q\}$.\\
Depth 2: &point $p_1=l\cap r_1$ and point $p_2=l\cap v_2$.\\
Depth 3: &horizontal line $h_2$ passing through $p_1$.\\
Depth 4: &point $p_3=h_2\cap v_1$.\\
Depth 5: &line $r_2$ passing through $\{p_2,p_3\}$.\\
Depth 6: &point $p_4=r_2\cap h_1$.\\
Depth 7: &line $l'$ passing through $\{a,p_4\}$.
\end{tabular}

\begin{figure}
\begin{center}
\includegraphics[width=0.64\textwidth]{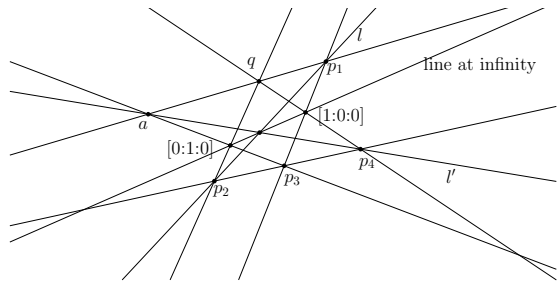}
\caption{How to construct a parallel line through one point}
\label{fig:papus_paralela_proyect}
\end{center}
\end{figure}

In the algebraic case, if the input elements $a,l,q$ are generic, then the
construction yields a realization of the hypothesis of Pappus Theorem with one
of the lines being the line at infinity and two of the points are the points
at infinity with projective coordinates $[0:1:0]$ and $[1:0:0]$, see
Figure~\ref{fig:papus_paralela_proyect}. Pappus theorem implies that the lines
$l$, $l'$ intersects at the line at infinity. Thus, $l'$ is the parallel to $l$
passing through $a$. The same approach work if we replace $l$ (a generic line) by a line passing through the affine origin of coordinates (curve of support $\{(1,0),(0,1)\}$) and $a$. We will use this construction as an auxiliary for the following:

Take as input points $a,b,c,q$, let $o=(0,0)$ be the origin of coordinates in
the affine plane $\mathbb{K}^2$, a line through a point $p$ and $o$ is just the
curve through $p$ of support $\{(1,0),(0,1)\}$. Consider the following
construction:\\
\begin{tabular}{ll}
Depth 1: &$l_1=\overline{oa}$, $l_2=\overline{ob}$, $l_3=\overline{oc}$\\
Depth 2-8: &$l_4=Parallel(a,l_2,q)$, $l_5=Parallel(b,l_1,q)$\\
Depth 9: &$d=l_4\cap l_5$\\
Depth 10: &$l_6=\overline{od}$\\
Depth 11-17: &$l_7=Parallel(d,l_3,q)$, $l_8=Parallel(c,l_6,q)$\\
Depth 18: &$z=l_7\cap l_8$\\
Depth 19: &$l_9=\overline{az}$
\end{tabular}

In the affine plane, we have constructed the parallelograms $oadb$ and $odzc$.
Hence, if $a=(a_1,a_2)$, $b=(b_1,b_2)$ and $c=(c_1,c_2)$, then $d=(a_1 +b_1, a_2
+b_2)$ and $z=(a_1 +b_1 +c_1, a_2 +b_2 +c_2)$. Notice that this construction if
far from being an admissible one.

Take the following tropical input elements of this construction, $a=(0,0)$,
$b=(-1,-1)$, $c=(-2,-2)$ and $q=(2,-1)$. For this input, we have that $z=(0,0)$
and $l_9=``0x+0y+0"$. The constructible set $\mathfrak{S}$ associated to this
input is the empty set. Lifts of the input elements are
\[\widetilde{a}=(\alpha_1+\ldots, \alpha_2+\ldots),
\widetilde{b}=(\beta_1t+\ldots, \beta_2 t
+\ldots),\]
\[\widetilde{c}=(\gamma_1 t^2 +\ldots, \gamma_2 t^2 +\ldots), \widetilde{q}=(
\eta_1 t^{-2}
+\ldots, \eta_2 t +\ldots )\]

The algebraic computations of $\widetilde{z}$ leads to the point
\[\widetilde{z}=(\alpha_1+\ldots, \alpha_2+\ldots).\]
That is, the principal term of $\widetilde{a}$ and $\widetilde{z}$ are the same.
So, we cannot compute the algebraic line $\widetilde{l}_9$ neither we cannot
deduce if the generic lift of the input will work or if there will be a lift at
all. However, it can be checked that the set $\mathfrak{S}_z$ associated to the
subconstruction that defines $z$ is nonempty and dense $\{\beta_2 -\eta_2\neq 0,
\alpha_2 \beta_1 -\alpha_1 \beta_2\neq 0, -\alpha_1 \gamma_2 +\gamma_1
\alpha_2\neq 0\}\cap (k^*)^8$.

In fact, for this construction and this tropical realization, the generic lift
works and it is compatible with tropicalization. To explain this, we know that
$\widetilde{z}=\widetilde{a}+\widetilde{b}+\widetilde{c}$. If
$\widetilde{a}=(\widetilde{a}'_1, \widetilde{a}'_2),
\widetilde{b}=(\widetilde{b}'_1t, \widetilde{b}_2't),
\widetilde{c}=(\widetilde{c}'_1t^2, \widetilde{c}_2' t^2),
\widetilde{q}=(\widetilde{q}'_1 t^{-2}, \widetilde{q}'_2 t)$, where
$\widetilde{a}_i', \widetilde{b}_i', \widetilde{c}_i', \widetilde{q}_i'$ are
elements of valuation zero. Then $\widetilde{z}=(\widetilde{a}'_1
+\widetilde{b}'_1t+\widetilde{c}'_1t^2, \widetilde{a}'_2
+\widetilde{b}'_2t+\widetilde{c}'_2t^2)$ and $\widetilde{l}_9=(\widetilde{b}'_2
t+\widetilde{c}'_2t^2)x+(-\widetilde{b}'_1 t -\widetilde{c}_1 t^2)y +
(\widetilde{a}'_2 \widetilde{b}'_1 -\widetilde{a}'_1 \widetilde{b}'_2)t
+(\widetilde{a}'_2 \widetilde{c}'_1 -\widetilde{a}'_1 \widetilde{c}'_2) t^2=0$.
If $\alpha_2\beta_1 -\alpha_1\beta_2\neq 0$ then
$T(\widetilde{l}_9)=``(-1)x+(-1)y+(-1)"=``0x+0y+0"=l_9$.

As a negative example, take the same construction but we take as input element
$b=(-1,-2)$, then we will arrive to the same situation of undecidability as
above, the set $\mathfrak{S}$ is again empty. If we take as before generic lifts
of the input elements, but this time $\widetilde{b}=(\widetilde{b}'_1t,
\widetilde{b}'_2t^2)$. Now, $\widetilde{z}=(\widetilde{a}'_1+\widetilde{b}'_1 t
+\widetilde{c}'_1t^2, \widetilde{a}'_2+(\widetilde{b}'_2 +\widetilde{c}'_2)t^2)$
and $\widetilde{l}_9=(\widetilde{b}'_2 +\widetilde{c}'_2) t x +(-\widetilde{b}'_1
-\widetilde{c}'_1 t)y+ \widetilde{a}'_2 \widetilde{b}'_1 + (\widetilde{a}'_2
\widetilde{c}'_1 - \widetilde{a}'_1 \widetilde{b}'_2 - \widetilde{a}'_1
\widetilde{c}'_2 )t$. Then $T(\widetilde{l}_9)=``(-1)x+0y+r"$, where $r\geq 0$.
So it never tropicalizes correctly.
\end{exmp}

\section{Notion of Constructible Theorem}\label{section_theorem}

Many classical theorems in Projective Geometry deal with
properties of configurations of points and curves. Thus, we can use the
relationship between the algebraic and tropical configurations in order to
transfer a Theorem from Classical Geometry to Tropical Geometry. So, we need a
notion of \emph{``Theorem''} is terms of configurations. We propose the
following notion.

\begin{defn}\label{teorema_construible}
A \emph{constructible incidence statement} is a triple $(G,H,x)$ such that $G$
is an incidence structure, $H$ is a geometric construction, called the
\emph{hypothesis}, such that, considered as an incidence configuration, $H$ is a
full substructure of $G$, $H \subseteq G$. Moreover,
\[\{\mathfrak{p}_G \cup \mathfrak{B}_G \} \setminus \{\mathfrak{p}_H \cup
\mathfrak{B}_H \} = \{x\},\]
there is only one vertex $x$ of $G$ which is not a vertex of $H$, this is called
the \emph{thesis node}.

Let $H_0$ be the set of input elements of $H$ as a construction. Let $\mathbb{K}$ be
an algebraically closed field. The incidence statement \emph{holds} in $\mathbb{K}$ or
it is a \emph{constructible incidence theorem} over $\mathbb{K}$ if it holds for
the generic realization of $H_0$. That is, if there is a non empty open set $L$
defined in the support space of $H_0$, $L\subseteq S_{H_0}$ such that:
\begin{itemize}
\item For every $\widetilde{h}\in L$, the construction $H$ is well defined.
\item If $\widetilde{p}\in R_H$ is the realization of $H$ constructed from $\widetilde{h}$, then
there is an element $\widetilde{x}$ such that $(\widetilde{p},\widetilde{x})$ is a realization of $G$.
\end{itemize}

In the tropical context, the construction $H$ is always well defined. Every
realization $h$ of the input of $H$ defines a realization $p$ of $H$ by the
construction. So, a constructible statement \emph{holds} in the tropical plane
or it is a \emph{tropical constructible incidence theorem} if, for each
realization $p$ of $H$ obtained by the construction, there is a tropical element
$x$ such that $(p,x)$ is a tropical realization of $G$.
\end{defn}

\begin{exmp}
There are many straightforward theorems that fit in this definition. For
example, let $H_0=\{p_1,p_2,l_1\}$, where $p_1,p_2$ are points and $l_1$ is a
line. Let $\mathfrak{C}$ be the construction consisting in computing the line
$l_2$ through $p_1$ and $p_2$. Let $x$ be the thesis node representing a point
and impose the conditions that $x$ belongs to both lines $l_1$ and $l_2$. The
vertices of $G$ are $\{p_1,p_2,l_1,l_2\}$. The edges (incidence conditions) of
$G$ are those of $H$, $\{(p_1,l_2), (p_2,l_2)\}$ plus the edges connecting the
thesis node $\{(x,l_1), (x,l_2)\}$. This statement only asserts that $l_1$,
$l_2$ have a common point. So it holds in every field $\mathbb{K}$ and also in the
tropical plane $\mathbb{T}^2$.
\end{exmp}

Of course, this notion is interesting if the thesis node $x$ and the elements
linked to it $h_1,\ldots, h_n$ form an incidence structure $G_0$ that is not
realizable whenever the elements $h_1,\ldots,h_n$ are generic. For instance, the
case where $x$ is a line containing three points $h_1$, $h_2$ and $h_3$. Now we
prove a transfer result for constructible incidence theorems.
\begin{thm}\label{transfer_theorem}
Let $\mathcal{Z}=(G,H,x)$ be a constructible incidence statement. Suppose that
the construction $H$ is admissible. If $\mathcal{Z}$ holds in a concrete
algebraically closed field $\mathbb{K}$, then it holds for every tropical plane
$\mathbb{T}^2$.
\end{thm}
\begin{proof}
First, suppose that $\mathbb{T}$ is the value group of the algebraically closed
field $\mathbb{K}$ such that $\mathcal{Z}$ holds. Let $h$ be a tropical realization
of the input elements of the hypothesis $H$. Let $p$ be the tropical realization
of $H$ constructed from $h$. As $H$ is an admissible construction, by
Theorem~\ref{admisible_es_levantable}, the set $\mathfrak{S}$ defined in
$(k^*)^N$
associated to $h$ contains a non empty open set. It follows that there is always
a lift $\widetilde{h}$ of $h$ belonging to $L$ and such that its principal coefficients
belong to the set $\mathfrak{S}$. Then, we can lift $p$ to an algebraic
realization $\widetilde{p}$ of $H$ constructed from $\widetilde{h}$. As $\mathcal{Z}$ holds in
$\mathbb{K}$, there is an element $\widetilde{x}$ such that $(\widetilde{p},\widetilde{x})$ is a realization
of $G$. It follows that its projection $(p,x)$ is a tropical realization of $G$
and $\mathcal{Z}$ holds in $\mathbb{T}$.

For the general case, the set $L$ of good input elements of $H$ is definable in
the first order language of the prime field of $\mathbb{K}$. So, if the theorem
holds in an algebraically closed field, it holds over any algebraically closed
field of the same characteristic \cite{Robinson}. In particular, fixed a
tropical semifield $\mathbb{T}$, there is an algebraically closed valued field
$\mathbb{L}$ of the same characteristic as $\mathbb{K}$ and whose valuation group is
$\mathbb{T}$. Thus, if $\mathcal{Z}$ holds in $\mathbb{K}$, then it also holds in
$\mathbb{L}$ and hence, it holds in $\mathbb{T}$.
\end{proof}

\subsection{Examples of Theorems}
Some examples of constructible incidence theorems are shown. They are all classic, but they are rewritten as constructible incidence theorems. There is an additional problem when expressing the theorems this way. Usually, it is not enough to provide a naive construction of the hypothesis, because it is very likely that the resulting construction is not admissible and Theorem~\ref{transfer_theorem} does not apply. So, the presentation of the theorems might seem strange at first sight.

\subsubsection{Fano Plane Configuration Theorem}
This first example shows the dependence of the characteristic of the field
$\mathbb{K}$ in order to derive the validity of a constructible incidence theorem in
the tropical context. The classical Theorem deals with the configuration of
points and lines in Fano plane, the projective plane over the field $\mathbb{F}_2$.
\begin{figure}
\begin{center}
\includegraphics[width=0.25\textwidth]{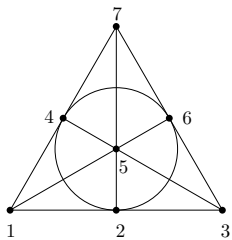}
\caption{The configuration of Fano plane}
\label{fig:Fano_plane}
\end{center}
\end{figure}
The configuration of Fano plane consists in 7 lines and 7 points as represented
in Figure~\ref{fig:Fano_plane}. This configuration cannot be realized over a
plane of characteristic zero. In a field of
characteristic $2$, if seven points $1, 2, 3, 4, 5, 6, 7$ verifies that the
triples $(1,2,3)$, $(1,4,7)$, $(3,6,7)$, $(1,5,6)$, $(2,5,7)$, $(1,4,7)$ are
collinear, then the points $(2,4,6)$ are also collinear. This Theorem holds in a
field $\mathbb{K}$ if and only if the field is of characteristic $2$. About the
tropicalization of this Theorem, it was proved to hold in $\mathbb{T}^2$ by M.
Vigeland using specific techniques \cite{vigeland-fano}. See also
\cite{MR2178322} for an application of this configuration to the comparison of
different notions of the tropical rank of a tropical matrix.
\begin{thm}[Fano plane configuration
Theorem]\label{T_FANO}\mbox{}\\
Construction of the hypothesis $H$:\\
\begin{tabular}{ll}
Input: & points $1,2,3,4,5$.\\
Depth 1: & lines $a=\overline{13}$, $b=\overline{15}$, $c=\overline{17}$,
$d=\overline{35}$, $e=\overline{37}$, $f=\overline{57}$.\\
Depth 2: & points $2=a\cap f$, $4= c\cap d$, $6= b \cap e$.\\
Thesis node: & $l$
\end{tabular}\\
Thesis: points $2,4,6$ belong to $l$.\\
The construction of the hypothesis is admissible, so we can derive that the
theorem holds in the tropical plane. In brief, this Theorem proves that, if we
start with any set of points $1$, $3$, $5$, $7$ in which even we may allow
repetitions and we perform the construction steps above, then three new points
$2,4,6$ will be obtained, and these three new points will necessarily lie on a
common tropical line $l$.
\end{thm}

\subsubsection{Pappus Theorem}
This classical theorem was studied from a tropical perspective in \cite{RGST}.
There, the authors showed that a direct translation of the usual hypothesis of
the theorem does not imply the thesis in the tropical context. On the other
hand, they proposed a constructive version of this Theorem. We proved this
constructive version of this Theorem in \cite{Pappus-trop} using a precursor
technique of our construction method.

\begin{thm}[Pappus Theorem]\label{T_PAPPUS}\mbox{}\\
Construction of the hypothesis $H$:\\
\begin{tabular}{ll}
Input: & points $1,2,3,4,5$.\\
Depth 1: & lines $a=\overline{14}$, $b=\overline{24}$, $c=\overline{34}$,
$a'=\overline{15}$, $b'=\overline{25}$, $c'=\overline{35}$.\\
Depth 2: & points $6=b\cap c'$, $7= a'\cap c$, $8= a \cap b'$.\\
Depth 3: & lines $a''=\overline{16}$, $b''=\overline{27}$,
$c''=\overline{38}$.\\
Thesis node: &point $p$
\end{tabular}\\
Thesis: lines $a''$, $b''$, $c''$ pass through $p$.
\end{thm}

\subsubsection{Converse Pascal Theorem}
Let $A$, $B$, $C$, $A'$, $B'$, $C'$ be six points in the plane, let
$P=\overline{AB'}\cap \overline{A'B}$, $Q=\overline{BC'}\cap \overline{B'C}$,
$R=\overline{AC'}\cap \overline{A'C}$. Converse Pascal Theorem proves that if
$P,Q$ and $R$ are collinear, then $A,B,C,A',B',C'$ belong to a conic. The
dimension of the space of realizations of a Pascal configuration is 11: 5
degrees of freedom comes from the conic and the points $A,B,C,A',B',C'$
belonging to the conic adds one degree of freedom each. If we want to define a
constructible theorem such that the thesis node is the conic, then the algebraic
elements of the construction of the hypothesis can only be points and lines. By
the nature of the steps of a construction, any construction that only uses
points and lines will provide configurations whose realization space has even
dimension (as it equals the dimension of the support space of the input
elements). It follows that the dimension of the support space of any potential
construction of a Pascal configuration $H$ is even. So, we cannot obtain such a
construction for this theorem. However, we can define a bigger construction such
that it contains Pascal configuration as a substructure. Namely, we can add
three arbitrary points points $X_1$, $X_2$, $X_3$ belonging to $\overline{AB'}$,
$\overline{BC'}$, $\overline{CA'}$ respectively, see
Figure~\ref{fig:Pascal_conv_teo}. Hence our configuration $G$ is Pascal
configuration with three additional marked points $X_1$, $X_2$, $X_3$. Its
dimension is now $14$. This is a example of how an additional step
\emph{``choose a line through $A$"} in a construction can be modeled by adding
the additional free point $X_1$ and then defining the line $\overline{AX_1}$.
\begin{thm}[Converse Pascal
Theorem]\label{pascal_converse}\mbox{}\\
\begin{figure}
\begin{center}
\includegraphics[width=0.5\linewidth]{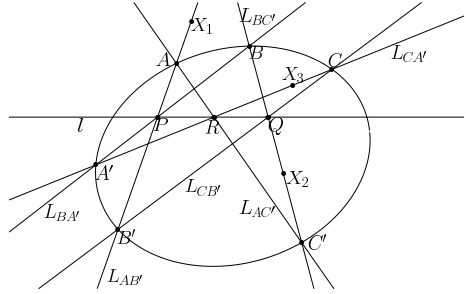}
\caption{Converse Pascal Theorem}\label{fig:Pascal_conv_teo}
\end{center}
\end{figure}
Construction of the hypothesis $H$:\\
\begin{tabular}{ll}
Input: & points $A,B,C, X_1,X_2,X_3$, line $l$.\\
Depth 1: & lines $L_{AB'}=\overline{AX_1}$, $L_{BC'}=\overline{BX_2}$,
$L_{CA'}=\overline{CX_3}$.\\
Depth 2: & points $P=L_{AB'}\cap l$, $Q=L_{BC'} \cap l$, $R=L_{CA'} \cap
l$.\\
Depth 3: & lines $L_{AC'}=\overline{AR}$, $L_{BA'}= \overline{BP} $,
$L_{CB'}= \overline{CQ}$.\\
Depth 4: & points $A'=L_{CA'}\cap L_{BA'}, B'=L_{AB'}\cap L_{CB'},$\\
& $C'=L_{AC'}\cap L_{BC'}$.\\
Thesis node: & conic $R$.
\end{tabular}\\
Thesis: points $A, B, C, A', B', C'$ belong to conic $R$.
\end{thm}
\subsubsection{Chasles Theorem}
Chasles Theorem \cite{MR1376653} states that if $\{q_1,\ldots, q_9\}$ are
the intersection points of two cubics, then any cubic passing through
$\{q_1,\ldots, q_8\}$ also passes through $q_9$. This implies that given another
free point $q_{0}$, there is always a cubic through $\{q_0,q_1,\ldots,q_9\}$.
This version can be easily translated to the tropical context.
\begin{thm}[Chasles Theorem]\label{T_CHASLES}\mbox{}\\
Construction of the hypothesis $H$:\\
\begin{tabular}{ll}
Input: & cubics $C_1$, $C_2$, point $q_0$.\\
Depth 1: & points $\{q_1, \ldots, q_9\}=C_1 \cap C_2$.\\
Thesis node: & cubic $R$.
\end{tabular}\\
Thesis: points $\{q_0, q_1, \ldots, q_9\}$ belong to cubic $R$.
\end{thm}

It is not true that every cubic passing through eight of the intersection points
passes through the ninth. See Figure~\ref{fig:Chasles_counterexample}. Let
$f=``0 +1x +1y
+1x^2 +3xy +1y^2 +0x^3 +1x^2y +1xy^2 +0y^3$, $g=19 +14x +20xy +24y +7x^2 +12x^2y
+23xy^2 +28y^2 +0x^3 +31y^3"$,
\[\begin{matrix}
f\cap_{st} g=&\{(-1, -3), &(0, -3), &(1, -3),\phantom{\}}\\&\phantom{\{}(-1,
-4), &(0, -4), &(1, -4),\phantom{\}}\\&\phantom{\{}(-1,-5), &(0, -5), &(1,
-5)\phantom{,}\}
\end{matrix}
\]
Take $h=``0 +1x +5y +\frac{11}{2}xy +1x^2 +9y^2 +5x^2y +9xy^2 +0x^3 +12y^3"$.
This is a cubic passing through $8$ of the stable intersection points of $f$ and
$g$ but not through the ninth.

An alternative to the Chasles Theorem that also holds in the tropical plane is
the following. Take as $8+n$ points $\{q_1,\ldots,q_8\}$, $\{x_1, \ldots,
x_n\}$, $n\geq 3$. All the steps are computing the cubic $C_i$ passing through
$\{q_1,\ldots, q_8, x_i\}$, $1\leq i\leq n$. The thesis node is a point $x$ and
the thesis is that $x$ belongs to $C_i$, $1\leq i\leq n$. The difference with
the previous version of Chasles theorem is that, by construction, the eight
points $\{q_1,\ldots, q_8\}$ are always in general position in every cubic
$C_i$. In our example, the points are not in general position neither in
$\mathcal{T}(f)$ nor $\mathcal{T}(g)$.

An immediate generalization of Chasles Theorem is the following.
\subsubsection{Cayley-Bacharach Theorem}
The generalization of Chasles Theorem (cf \cite{MR1376653}) we discuss here is
the following:
let $C_1$, $C_2$ be plane curves of degrees d and e respectively, intersecting
in $de$ distinct points $Q=\{p_1 ,\ldots, p_{de}\}$. If $C$ is any plane curve
of degree $d + e - 3$ containing all but one point of $Q$, then $C$ contains
every point of $Q$. The second version of Chasles Theorem given does not fit
well to this theorem, but the generalization of the first version of Chasles
Theorem is immediate, note that a curve of $d+e-3$ is determined by
$\frac{d^2+e^2-3e-3d}{2}$ points:

\begin{figure}
\begin{center}
\includegraphics[width=0.4\linewidth]{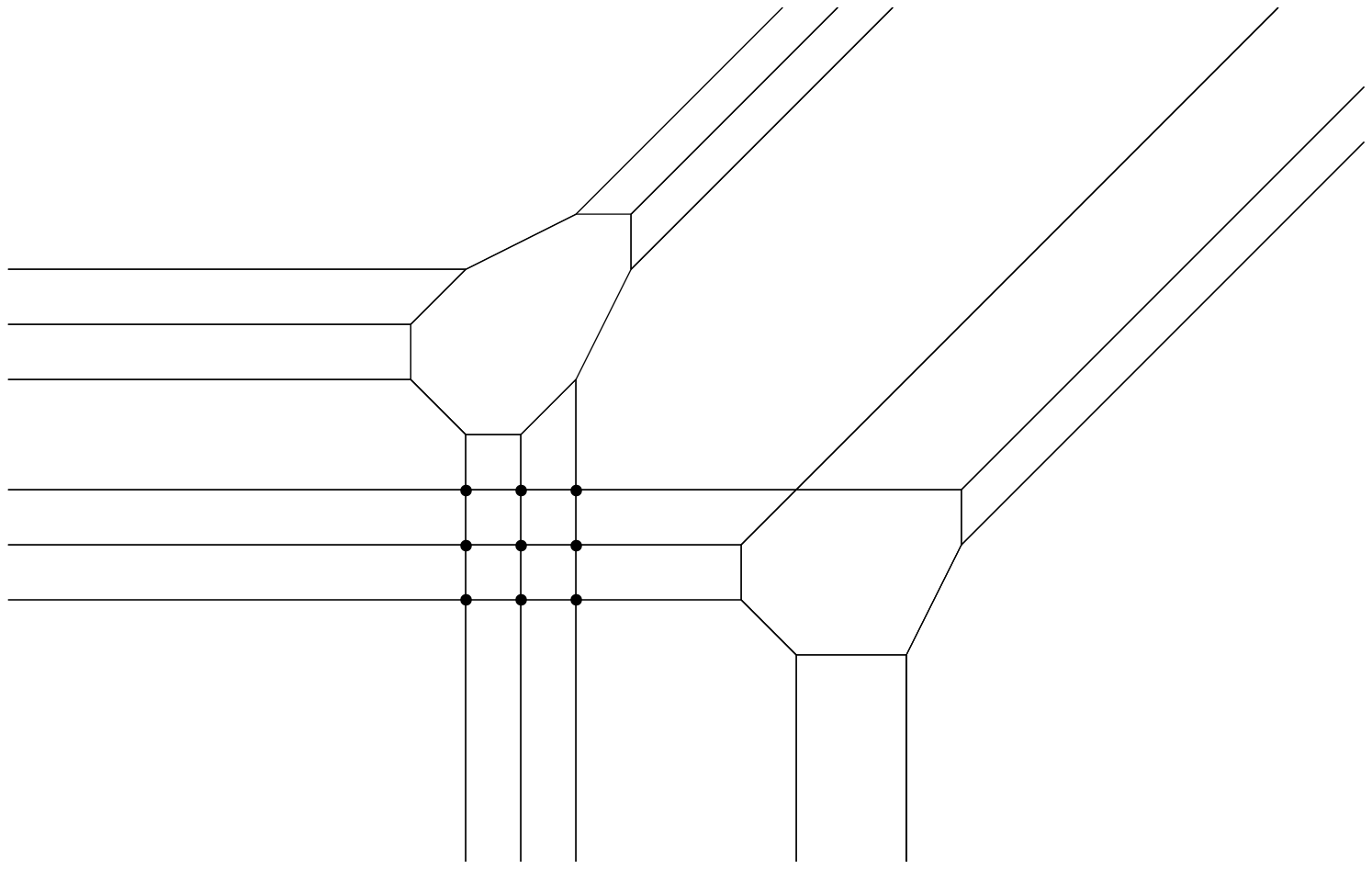}
\includegraphics[width=0.3\linewidth]{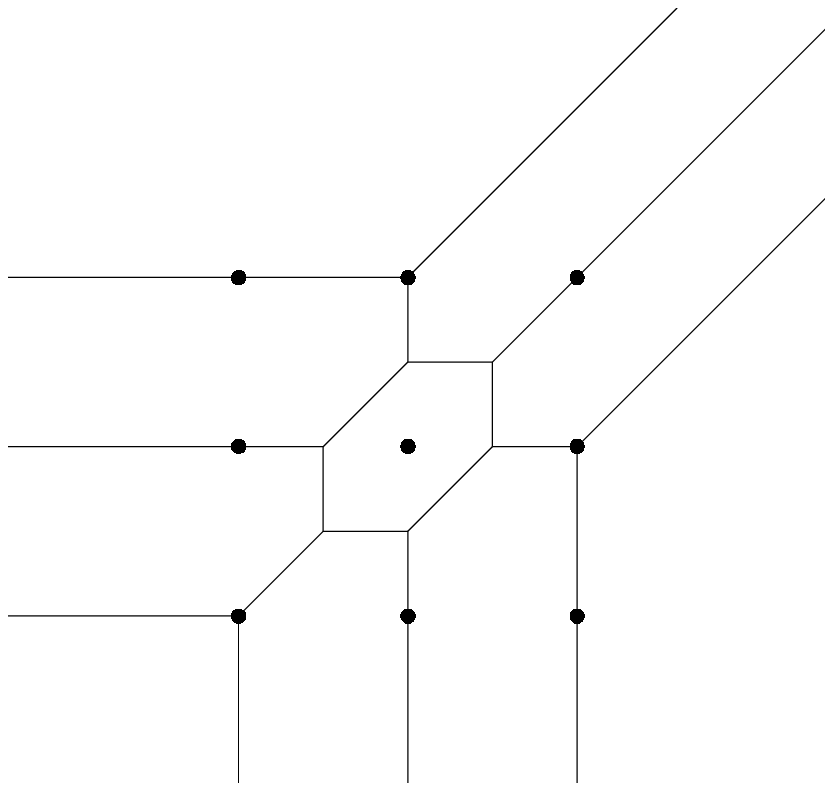}
\caption{A cubic through 8 but not 9 intersection points of two other cubics}
\label{fig:Chasles_counterexample}
\end{center}
\end{figure}
Let $d,e\geq 3$ natural numbers, $l=1+\frac{d^2+e^2-3e-3d}{2}$
\begin{thm}[Cayley-Bacharach Theorem]
\label{T_CB}\mbox{}\\
Construction of the hypothesis $H$:\\
\begin{tabular}{ll}
Input: & degree $d$ curve $C_1$, degree $e$ curve $C_2$, points
$p_1,\ldots, p_l$.\\
Depth 1: & points $\{q_1, \ldots, q_{de}\}=C_1 \cap C_2$.\\
Thesis node: & curve $R$ of degree $d+e-3$.
\end{tabular}\\
Thesis: points $\{q_1, \ldots, q_{de}\} \cup \{p_1,\ldots, p_l\}$ belong to
curve $R$.
\end{thm}
\subsubsection{Weak Pascal Theorem}
This Theorem is not in the context of Theorem~\ref{transfer_theorem} because the
construction involved is not admissible. Nevertheless, for some tropical
realization of the hypothesis, we will be in the context of
Theorem~\ref{Construccion_casi_admisible}. So this Theorem does not hold for
every tropical input, we have to add conditions in the tropical realization.

\begin{thm}[Weak Pascal Theorem]\label{weak_pascal}\mbox{}\\
Consider the following construction:\\
\begin{tabular}{ll}
Input: & conic $Z$, lines $L_1$, $L_2$, $L_3$.\\
Depth 1: & points $\{A,B'\}=R \cap L_1$, $\{B,C'\}=R\cap L_2$, $\{C, A'\}=R
\cap L_3$.\\
Depth 2: & lines $L_4=\overline{AC'}$, $L_5=\overline{BA'}$,
$L_6=\overline{CB'}$.\\
Depth 3: & points $P=L_1 \cap L_5$, $Q= L_2\cap L_6$, $R=L_3\cap L_4$.
\end{tabular}\\
If a tropical instance of this construction is such that each set of points
$\{A, C'\}$, $\{B, A'\}$ and $\{C, B'\}$ is in generic position with respect to
$Z$, then there is a line $L$ (thesis node) that contains the points $P$, $Q$ and $R$.
\end{thm}
\begin{proof}
This construction, in the algebraic context, provides instances of Pascal
theorem. Hence, if the input is generic, then points $\widetilde{P}$, $\widetilde{Q}$,
$\widetilde{R}$ are collinear. But this construction is not admissible, so
Theorem~\ref{transfer_theorem} does not apply. Nevertheless, this construction
is in the
context of Theorem~\ref{Construccion_casi_admisible}. The minimal multiples
paths are $Z \rightrightarrows L_4$, $Z \rightrightarrows L_5$ and $Z
\rightrightarrows L_6$. By Theorem~\ref{Construccion_casi_admisible}, if each
one of these three sets is in general position with respect to $R$, then this
tropical instance can be lifted to the a generic instance in the
algebraic framework. As Pascal Theorem holds in $\mathbb{K}$. $\widetilde{P}$, $\widetilde{Q}$ and
$\widetilde{R}$ are collinear. So $P$, $Q$ and $R$ will be collinear.
\end{proof}

\begin{exmp}
Let $Z=``3y+5+3y^2+0x^2+4x+0xy"$ $L_1=``1y+0x+0"$ $L_2=``0y+0x+2"$
$L_3=``(9/2)y+0x+3"$, then $A=(3,2)$, $B'=(1,0)$, $B=C'=(2,3/2)$, $C=(1,-3/2)$,
$A'=(4,-1/2)$, $L_4=``3y+2x+(9/2)"$, $L_5=(3/2)x+4y+(11/2)$, $L_6=0x+1y+1$,
$P=(5/2,3/2)$, $Q=(2,1)$, $R=(5/2,-3/2)$. The points $P$, $Q$ and $R$ are not
collinear, in this example, the set $\{C,B'\}$ is not in generic position in
$Z$. 

However, for these input elements, the election of the points in the depth 1
steps is arbitrary. If we now take $A=(1,0)$, $B'=(3,2)$, $B=C'=(2,3/2)$,
$C=(4,-1/2)$ and $A'=(1,-3/2)$, now $L_4=``2y+(3/2)x+(5/2)"$,
$L_5=``2y+(3/2)x+(5/2)"$, $L_6=``4y+2x+6"$, $P=(1,0)$, $Q=(2,2)$, $R=(1,-3/2)$.
In this case, the three sets of points are in generic position in $Z$, it can be
checked that the three points belong to the tropical line of equation
$L=``2x+2y+3"$.
\end{exmp}

\medskip
Luis Felipe Tabera Alonso\\
Dept. \`Algebra i Geometria, Facultat de Matem\`atiques, Universitat de Barcelona\\
Barcelona, Spain\\
e-mail: tabera@ub.edu

\end{document}